\documentclass[11pt]{article}

\usepackage {ucs}
\usepackage [utf8x]{inputenc}
\usepackage {graphicx}

\usepackage {amsfonts,amsmath,amssymb}
\usepackage {graphics, graphicx}
\usepackage {multicol}						
\usepackage {multirow}  					
\usepackage {soul} 							
\usepackage {datetime} 						
\usepackage {parskip}						
\usepackage {setspace}						
\usepackage{mathtools}						
\usepackage{pgfplots}
\usepackage{algorithmic}				
\usepackage{booktabs}

\usepackage{algorithm}
\usepackage{amsthm}

\newtheorem{definition}{Definition}

\usetikzlibrary{positioning}
\usetikzlibrary{arrows}
\usetikzlibrary{shapes}
\pgfplotsset{compat=newest}
\usetikzlibrary{plotmarks}
\usetikzlibrary{arrows.meta}
\usetikzlibrary{matrix}
\usetikzlibrary{pgfplots.groupplots}
\pgfplotsset{plot coordinates/math parser=false}
\pgfplotsset{try min ticks=3}
\pgfplotsset{plot coordinates/math parser=false}
\usepgfplotslibrary{colormaps} 
\usetikzlibrary{pgfplots.colormaps}
\usepgfplotslibrary{colorbrewer} 
\usetikzlibrary{colorbrewer}
\newlength\figureheight
\newlength\figurewidth

\usepackage{filemod}						

\renewcommand{\(}{\left(}
\renewcommand{\)}{\right)}
\renewcommand{\vec}{\mathbf}

\newcommand{\pd}{\partial}

\newcommand {\N}{\mathbb{N}}
\newcommand {\R}{\mathbb{R}}

\newcommand {\lmatrix}{\left[\begin{matrix}}
\newcommand {\rmatrix}{\end{matrix}\right]}

\newcommand {\startenv} {\vskip 0.05em \begin{tabular}{||l}\parbox[t]{0.95\linewidth}}
\newcommand {\stopenv} {\end{tabular}\vskip 0.05em}

\newcommand{\change}[2]{{#2}}

\begin{document}
\title{An adaptive rectangular mesh administration and refinement technique with application in cancer invasion models}
\author{Niklas Kolbe\textsuperscript{1}
	\and
	Nikolaos Sfakianakis\textsuperscript{2}\footnote{Corresponding author}
	}

\date{
  \scriptsize
  $^1$Institute of Geometry and Practical Mathematics, RWTH Aachen University, Germany \smallskip \\
  $^2$School of Mathematics and Statistics, University of St Andrews, UK \smallskip \\
  \tt{kolbe@igpm.rwth-aachen.de}, \ \
  \tt{n.sfakianakis@st-andrews.ac.uk}}
        
\maketitle

\begin{abstract}
	We present an administration technique for the bookkeeping of adaptive mesh refinement on (hyper-)rectangular meshes. Our technique is a unified approach for h-refinement on 1-, 2- and 3D domains, which is easy to use and avoids traversing the connectivity graph of the ancestry of mesh cells. Due to the employed rectangular mesh structure, the identification of the siblings and the neighbouring cells is greatly simplified. The administration technique is particularly designed for smooth meshes, where the smoothness is dynamically used in the matrix operations. It has a small memory footprint that makes it affordable for a wide range of mesh resolutions over a large class of problems. We present three applications of this technique, one of which addresses h-refinement and its benefits in a 2D tumour growth and invasion problem.
	
\end{abstract}


\section{Introduction}

Time and again \textit{adaptive mesh refinement} (AMR) has been employed to improve the accuracy of numerical methods and to alleviate their computational burden. Typical fields of application of AMR include engineering, astrophysics, and fluid dynamics, where the associated techniques have become a vital component of the overall numerical investigation, see, e.g.,~\cite{Berger1984, Berger1989,Babuska1978,Verfurth1994, Teyssier2002, Puppo2011,tenaud2011tutorials, sempliceAdaptiveMeshRefinement2016}. In mathematical biology however, they have been scarcely applied yet. Some examples can be found in~\cite{DudleyWard2011, Botti2010, Wise2008, Frieboes.2006, Kol-Kat-Sfak-Hel-Luk.2014}. 

AMR, has become a standard for which particular and elaborate numerical methods have been developed. Still the mesh administration remains a complex process, which most scientific computing efforts try to avoid. The aim of the current paper is to propose a ``do-it-yourself'' recipe of AMR and mesh administration strategy that emphasizes its simplicity of implementation and use, while allowing for the handling of situations that cover most of the spectrum of scientific computing.

In this perspective, one of the aims of our work is to apply AMR to models that describe the \textit{first stage in cancer metastasis} ---and one of the \textit{hallmarks of cancer}--- the invasion of the \textit{extracellular matrix} (ECM)~\cite{hanahan2011hallmarks}. The models we consider in this work are typically \textit{advection-reaction-diffusion} (ARD) systems where the involved quantities are represented by their macroscopic densities. The dynamics that their solutions exhibit can be complex and their numerical treatment is challenging, and hence, AMR is seen as an important tool. A 1D application of  AMR---in particular \textit{h-refinement}---in a \textit{Finite Volume} (FV) method for such a system has been studied in our previous work \cite{Kol-Kat-Sfak-Hel-Luk.2014}. There, we demonstrated a significant improvement in accuracy and efficiency when AMR was employed. For 2D (or higher dimensional) cases however, a series of additional numerical difficulties arise. For example, the nature of these particular cancer invasion models promotes the use of FV methods over meshes with \textit{rectangular cells}, which when combined with h-refinement techniques lead to \textit{hanging nodes}. This complicates significantly the computation of the advection and, most notably, of the diffusion components of the system. In the literature, 2D \textit{convection-diffusion} problems have been previously solved by FV methods augmented by AMR techniques, see e.g.~\cite{MacKenzie.2012} and the references therein. 

Many technical issues related to the implementation of the AMR techniques arise as well. Some are typical: traversing the \textit{connectivity graph} between \textit{children}, \textit{parents}, and \textit{ancestor} mesh cells, as well as the identification of \textit{sibling} and \textit{neighbour cells} needed for the computation of the numerical fluxes and for the refinement/coarsening procedures. Other difficulties are less common e.g. particularities of the model that might pose additional conditions on the structure and the handling of the mesh. To mention but a few: the fast dynamics of the solutions and the diffusion of the model manifest themselves in the ``smoothness'' of the mesh\change{;}{, we refer e.g. to \cite{edelsbrunner_2001} for a thorough discussion on the topology of mesh generation. Furthermore,} the large number of system components and the richness of phenomena to be resolved call for different (sometimes multiple) refinement and coarsening criteria; the possibility of a blow-up, which is inherent in this type of systems, \change{}{see e.g. \cite{Espejo2009, Winkler2021}, } requires a dynamic adaptation of the highest refinement level.

To address the above numerical and technical issues, we develop our own AMR and mesh administration techniques. There are several reasons for that: first, we aim for simplicity, flexibility, and portability of our algorithms. Secondly,  we want complete control of all the stages of the numerical treatment of the models. This includes various components of the implementation: the numerical methods for the solution of the \textit{partial differential equation} (PDE) systems, the AMR techniques for the refinement/coarsening process, as well as the grid administration and \textit{bookkeeping algorithms} for the handling of the data on the mesh.

The overall effort is extensive; we have previously addressed the development of the numerical methods to solve the corresponding models, and the 1D h-refinement technique~\cite{Hel-Kol-Sfak.2015, Sfak-Kol-Hel-Luk.2016, Kol-Kat-Sfak-Hel-Luk.2014}. In the current paper, we focus primarily on the grid administration/handling, and secondarily on the presentation of a concrete application of  h-refinement to a 2D \textit{cancer invasion} system. The \textit{h-refinement} technique that we employ, makes use of cell bisection in 2, 4, or 8 equivalent rectangle cells in 1-, 2-, or 3D, respectively. This approach has been widely used in the application of AMR methods, see e.g.~\cite{Puppo2011, sempliceAdaptiveMeshRefinement2016}. 

Regarding the structure of the mesh data and their bookkeeping, there are several requirements that should be fulfilled by a modern mesh administration technique. We mention here only the following\change{}{, and refer to \cite{Dobkin1989,Blandford2005,Alumbaugh2005,Canino2011, Dyedov2015} for more details}: a) \textit{generality} with respect to the dimensions of the problems and the shape of the domain, b) \textit{efficiency} in the access times of the stored information during and after the reconstruction of the mesh, c) at least affordable, if not minimal, \textit{memory costs}, d) \textit{simplicity, portability}, and \textit{transparency} of the implementation, and e) \textit{extensibility} to parallelizations. 

The most common practice \change{in the literature} is to use \textit{pointer-based mesh data structures}. In this approach, the basic objects considered---e.g. vertices---are explicitly represented by their physical coordinates. The edges/faces/cells are defined by reference to the vertices using a cascade of pointers/handles\change{.}{, see e.g. \cite{Alumbaugh2005, Dendro,Kremer2012}.} In the particular case of a 2D  triangulation, these references are to the three vertices as well as the three neighbouring triangles of each triangle are stored. This is an intuitive and relatively easy to implement and use technique, but it has a significant memory footprint, and poses additional computational burden at every step of the method\change{;}{, cf. \cite{Blandford2005}.} When a discretisation cell is refined the references pointing to this cell as well as the other cells it points to, need to be adjusted accordingly---a process that can be complex especially in 3D volumes. For further information on the implementation and the applications of this technique, we refer to~\cite{Kirk2006, Beall.1997, Carrey.1988, Becker1981, Kremer2013, Dyedov2015} and the references therein. Alternative approaches have been devised with main aim to increase the efficiency of the methods, see e.g. the \textit{half-edge} and \textit{half-facet array-based mesh data structure} implementations\change{.}{, see e.g.~\cite{Alumbaugh2005, Canino2011, Dyedov2015}.} These techniques have been used primarily for the computer graphics representation of surface objects in 3D, see e.g.~\cite{CGAL2000, VTK2006, Dyedov2009}.

The approach we propose in the current paper is problem-specific and focuses on rectangular meshes in 1-, 2- and 3D. Unlike the pointer-based method, we store the full data structure of the mesh tree in the form of a matrix and refer to different mesh cells via the corresponding lines of the matrix. Every cell points to its children and parent cells. We do not store the physical coordinates, we instead store the indices of the mesh cells, i.e. the refinement level and an \textit{intra-level} identifier of every cell instead of the physical coordinates of its vertices. In more detail, for given minimum and maximum refinement levels, we pre-compile a matrix in which we store, for each cell of the discretisation tree, information on its parent and children cells. This along with an enforced \textit{grid regularity}, allows for the effortless and efficient refinement, coarsening, and neighbour identification processes. The memory footprint of the method we propose is reasonable, even for a large cascade of refinement levels in 2- and 3D experiments without using memory compression.

When comparing with the software library \texttt{p4est} \cite{p4estSA}, our approach exhibits both commonalities and differences. On the one hand, both approaches assume a tree based cell-hierarchy that branches-off from a root cell and reaches the leaf cells. Each cell points to its parent and its children cells. On the other hand, in our approach we characterize the computational cells in the mesh tree by a single integer (in one-, two-, and three-D), where in \texttt{p4est} every 3-dimensional mesh cell is characterized by 3 (real number) coordinates $x$, $y$, $z$. Moreover, our representation refers to the centre of the computational cell, whereas in \texttt{p4est} the lower left corner is represented. This introduces some ambiguity as no two different cells in a binary-, quad-, or octree share the same centre but they might share the same lower left corner (as, e.g., parent and some child cells). \change{}{ For more information on the structure and functionality of \texttt{p4est} we refer to \cite{p4estSA,p4est15}. We also refer to it's employment in the \texttt{deal.II} computational framework \cite{deal.II, deal.IIb}.} The comparison \change{}{of our proposed technique} with \texttt{alugrid},~\cite{alugrid} is less straight forward as the information stored in \texttt{alugrid} is \change{}{more} diverse and includes vertices, edges, faces, and more. \change{}{Still, we refer the interested reader to \cite{alugrid, alugrid2} for a thorough discussion of this method and it's employment within the \texttt{DUNE} computational framework \cite{dune2008a, dune2008b}.}

Overall, in this work we propose a mesh data structure, which is easy to understand, implement, and use, and it can handle mesh resolutions that are sufficient for a wide range of academic investigations in 1-, 2- and 3D. We exhibit this flexibility by presenting three different applications/experiments while giving particular emphasis on a cancer invasion problem.

The rest of the paper is structured as follows: in Section~\ref{sec:mesh:admin} we present the mesh administration technique, address basic operations, such as the refinement/coarsening and the identification of neighbours and siblings, and discuss the memory usage and computational costs. In Section~\ref{sec:num:exper} we present three numerical experiments: a generic one exhibiting the properties of the mesh administration technique, as well as applications in a Euler system and a cancer invasion model. We discuss the basic ingredients of the numerical methods and the AMR technique we use.

\section{Mesh data structure and handling}\label{sec:mesh:admin}
A significant burden in the development of AMR techniques is the \textit{administration} of the discretisation meshes, in particular when dealing with time dependent 2- or 3D problems. We propose in this section a particularly easy procedure to \textit{represent}, \textit{store}, \textit{handle}, and in general \textit{bookkeep} \textit{dyadic rectangular meshes} in a unified way for 1-, 2-, or even 3D domains.

Main characteristic of the proposed mesh administration technique is that we encode the information of the dyadic tree inside an easy-to-use matrix. We do it in such a way that we access with ease the \textit{parent}, \textit{children}, and most notably the \textit{siblings} and \textit{neighbouring} cells of any given cell. In effect we simplify greatly the \textit{refinement} and, most importantly, the \textit{coarsening} steps of the mesh.

This part of the paper is structured as follows: in Section~\ref{sec:defs} we present the main notations and definitions, in Section~\ref{sec:cell:repr} we describe the way we store the information of the computational cells and the dyadic meshes, and comment on the memory usage of the method. In Section~\ref{sec:mat:op} we present the operations needed to identify the sibling and neighbour cells, in Section~\ref{sec:ref.coars} we elaborate on the refinement and coarsening procedures, and in Section~\ref{sec:projections} we discuss the handling of data on the mesh.

\subsection{Basic definitions and notations}\label{sec:defs}
The proposed technique can be employed on general hyperrectangles in 1-, 2- , or 3D domains. For the sake of simplicity we will restrict the presentation to the discretisation of the domain $\Omega=[0,1]^d$ for $d=1, \, 2$.

We denote by $G_l^d$, $l\in \N$ a \textit{uniform partition} of $\Omega$ with cardinality
\begin{equation}\label{eq:card}
	|G_l^d|=2^{ld}.
\end{equation}
Here, $l$ is the \textit{refinement level} of $G_l^d$ assumed to be bounded by  $l_\text{min}\leq l \leq l_\text{max}$ with $l_\text{min},l_\text{max}\in~ \N$. The elements of $G_l^d$ are called \textit{cells} and they are either intervals in $d=1$, or squares in $d=2$.   

For every cell $C \in \bigcup_{l=l_\text{min}}^{l_\text{max}} G_l^d$, we denote its refinement level, i.e. the index $l$ for which $C\in G_l^d$, by $L(C)$, its centre by $M(C)\in \Omega$, and the occupied physical domain by $D(C)\subset\Omega$. Accordingly we can write,
\begin{subequations}
	\begin{align}
		M(C) & \in \left\{ \sum_{i=1}^d \frac{2 k_i - 1}{2^{L(C)+1}} {\bf e_i},\quad 1 \leq k_i \leq 2^{L(C)}, \quad k_i \in \N \right\}, \label{eq:centr}         \\
		D(C) & = \left\{ M(C) +\sum_{i=1}^d \frac{\lambda_i}{2^{L(C)}} {\bf e_i},~-\frac12 < \lambda_i <\frac12, \quad \lambda_i\in\R \right\}. \label{eq:domn} 
	\end{align}
\end{subequations}
where ${\bf e_i}$, $i=1,2$ represents the unit vector of the corresponding axis. 

A cell $\tilde C$ is termed \textit{child cell} of the cell ${C} \in G_{l}^d$ if $\tilde{C} \in G_{l+1}^d$ and $M(\tilde{C}) \in D({C})$. Equivalently, the cell $C$ is called the \textit{parent cell} of $\tilde{{C}}$. Every parent cell has several (2,4, or 8, depending on the dimension $d$) children cells, and every child cell has a single parent cell. Children cells of the same parent cell are called \textit{siblings}. Geometrically, the children cells are obtained by \textit{bisection} of the parent cell.

The centre of a parent cell resides on the boundary of all its children cells, and the boundary of a parent cell is partly shared with the boundaries of it's children cells. 
The cell $\tilde{C} \in G_{\tilde l}^d$ is a \textit{descendant} of the cell $C \in G_{l}^d$, equivalent to saying that $C$ is an \textit{ancestor} of $\tilde{C}$, if there is a cascade of parent-children relations between $C$ and $\tilde C$.

Two cells that share a (part of their) boundary of co-dimension 1 but no part of their physical domain are called \textit{neighbours}. In 2D, this definition excludes cells that share a single point of their boundary. The common boundary between two neighbour cells is their \textit{interface}. 

For the rest of this work, we consider meshes that are subsets of $\bigcup_{l=l_\text{min}}^{l_\text{max}} G_l^d$, for particular $l_\text{min},\, l_\text{max}\in \N$. We also expect them to obey a particular smoothness relation discussed in the following definition.

\begin{definition}\textnormal{(Regular Structured Mesh).}
	A $d$-dimensional ($d=1,\,2$) mesh $G$ is called \textit{Regular Structured Mesh} (RSM) of minimum and maximum refinement levels $l_\text{min},l_\text{max}\in \N$, and of \textit{mesh regularity} $m_r\in \N$ if
	\begin{equation}\label{eq:grid}
		G \subset \bigcup_{l=l_\text{min}}^{l_\text{max}} G_l^d,
	\end{equation}
	and if 
	\begin{enumerate}
		\item For all ${\bf x}\in [0,1]^d$ either 	
		$\exists !\, C \in G$ such that ${\bf x} \in D(C)$ or $\exists C \in G$ such that ${\bf x} \in \pd {D(C)}$.
		
		\item For all neighbour cells $C_j, C_k \in G$, the \textit{mesh regularity condition} holds:
		\begin{equation}\label{eq:grid_reg}
			|L(C_j) - L(C_k)|\leq m_r.
		\end{equation}
	\end{enumerate}
\end{definition}

We can now introduce the basic operations of \textit{refinement} and \textit{coarsening}: for $G \subset \bigcup_{l=l_\text{min}}^{l_\text{max}} G_l^d$ a RSM of regularity $m_r$, and a cell $C\in G$, we set
\begin{itemize}
	\item[--] 
	\textit{Refinement} to be the process of replacing a cell $C$ in $G$ \textit{by all of its children} from the level $L(C)+1$. 
	\item[--] 
	\textit{Coarsening} to be the process of replacing a cell $C$ \textit{and all of its siblings} in $G$ by their parent cell from the level $L(C)-1$. 
\end{itemize}
As we expect the mesh $G$ to maintain the RSM properties after the refinement and coarsening, both operations are subject to additional constraints, i.e. the resulting cells have to be of refinement level $l$ such that $l_\text{min}\leq l\leq l_\text{max}$, and the resulting mesh should respect the mesh regularity condition \eqref{eq:grid_reg}. See Section \ref{sec:ref.coars}, for a detailed description of the refinement-coarsening procedures.


\subsection{Cell representation for nonuniform meshes}\label{sec:cell:repr}
The proposed mesh administration and book-keeping technique has the benefit of being a unified approach over 1-, 2-, and 3D domains, although for the sake of presentation we will discuss here only the 1- and 2D cases. 

In what follows we present the technique in its constituent components: the representation of cells and meshes, and how basic operations like refinement, coarsening or accessing siblings or neighbours are performed.
For ease of the presentation, we start with the 1D case. 

\subsubsection*{Cell representation in 1D}

\begin{table}[t]
	\begin{center}
		\begin{tabular}{cccc}\toprule
			level                        & number of  cells                    & centres                                                                                     & cell size                                     \\[0.5em]
			\midrule\\[-0.5em]
			0                            & 1                               & $\displaystyle \left\{\frac 1 2 \right\}$                                                   & 1                                         \\[1em]
			1                            & 2                               & $\displaystyle\left \{\frac 1 4, \frac 3 4 \right\}$                                        & $\displaystyle\frac 1 2$                  \\[1em]
			2                            & 4                               & $\displaystyle\left\{\frac 1 8, \frac 3 8 , \frac 5 8, \frac 7 8\right\}$                   & $\displaystyle\frac 1 4$                  \\[1em]
			$\vdots$                     & $\vdots$                        & $\vdots$                                                                                    & $\vdots$                                  \\
			$\displaystyle l_\text{max}$ & $\displaystyle2^{l_\text{max}}$ & $\displaystyle\left\{\frac{2k-1}{2^{l_\text{max}+1}}, ~k=1, \dots, 2^{l_\text{max}}\right\}$ & $\displaystyle\frac{1}{2^{l_\text{max}}}$\\ \bottomrule\\
		\end{tabular}
	\end{center}
	\caption{Mesh cells on uniform dyadic meshes that discretize [0,1] with refinement level $0$ through $l_\text{max}$. The refinement level $l$ and the intra-level index $k$ suffice for their complete characterization. }\label{tbl:cells}
\end{table}

We consider all possible mesh cells $C\in \bigcup_{l=l_\text{min}}^{l_\text{max}} G_l^1$ and identify their centres, relative positions, and sizes with respect to the levels of refinement in Table~\ref{tbl:cells}.				
We note at first that every cell $C$ can be uniquely defined by its refinement level $L(C)$ and its centre $M(C)$. The centre $M(C)$ in turn can be uniquely characterized by the refinement level $L(C)$ and an \textit{intra-level index} $k$ enumerating the cells in the current level, see also Table~\eqref{tbl:cells}.
				
Hence, the full sequence of dyadic meshes $\bigcup_{l=l_\text{min}}^{l_\text{max}} G_l$ can be described by the matrix
\begin{equation}\label{eq:C_1D}
  \mathcal C=\lmatrix
~&   \vec{k} &~\\
 ~&  \vec{l} &~ \\
 ~&  \vec{p} &~  \\
 ~&   \vec{d_l} &~ \\
~&   \vec{d_r} &~ \rmatrix^T,
\end{equation}
where the lines (i.e. the columns of $\mathcal{C}^T$) represent the cells $C\in \bigcup_{l=l_\text{min}}^{l_\text{max}} G_l^1$ and the line vectors $\vec{k}$, $\vec{l}$, and $\vec{p}$ include their intra-level index, their refinement level, and the line of $C$ in which the parent of the current cell is located. Furthermore $\vec{d_l}$, $\vec{d_r}$ are the lines where the children cells (left and right respectively) of the cell are stored in terms of line numbers of $\mathcal C$.  
Every cell included in the matrix $\mathcal C$ can be characterized uniquely by its refinement level $l$ and intra-level index $k$, or by the corresponding line of the matrix. 
				
Accordingly, the tree-example presented in Table \ref{tbl:cells} can be written in a  matrix formulation as 
\begin{equation}\label{eq_matrix_1D_example}
	\mathcal C=\left [\begin{array}{c |cc | cccc | cccccccc}
		1 & 1 & 2 & 1 & 2 & 3 & 4 & 1 & 2 & 3 & 4 & 5 & 6 & 7 & 8\\
		0 & 1 & 1 & 2 & 2  & 2 & 2 & 3 & 3 & 3 & 3 & 3 & 3 & 3 & 3\\
		- & 1 & 1 & 2 & 2 & 3 & 3 & 4 & 4 & 5 & 5 & 6 & 6 & 7 & 7\\
		2 & 4 & 6 & 8 & 10 & 12 & 14 & - & - & - & - & - & - & - & -\\
		3 & 5 & 7 & 9 & 11 & 13 & 15 & - & - & - & - & - & - & - & -\\
	\end{array} \right] ^T 
, 
\end{equation}
where the vertical lines separate cells of different refinement levels, with $l=3$ being the maximum in this case. The empty component in the first column of $\mathcal{C}^T$ implies that the corresponding first cell does not posses a parent cell within $\bigcup_{l=l_\text{min}}^{l_\text{max}} G_l^1$. Similarly, the empty components in the fourth and fifth line of $\mathcal{C}^T$ imply that the corresponding cells do not posses any children cells in $\bigcup_{l=l_\text{min}}^{l_\text{max}} G_l^1$.
				
\paragraph{Well posedness}
The centres of the cells of different levels cannot coincide since, otherwise, that would mean, for two cells with levels of refinement $l_1< l_2$ and intra cell indices $k_1, k_2$ that:
$$\frac{2k_1-1}{2^{l_1 +1}}=\frac{2k_2-1}{2^{l_2 + 1}}\quad{\Leftrightarrow}\quad \N\ni 2k_1-1=\frac{2k_2-1}{2^{l_2-l_1}}\not \in \N.$$
				
Since the size of every cell of level $l$ is $\frac{1}{2^{l}}$, cells of different levels will not intersect, unless one of them is a descendant of the other.

\paragraph{Size of $\mathcal C$ and memory usage in 1D}
As the refinement level $l$ is upper bounded by $l_\text{max}$, the size of the matrix $\mathcal C$ is finite. In particular, assuming a cascade of meshes starting from the coarse mesh $G_{l_\text{min}=5}^1$ of 32 cells up to the fine mesh $G_{l_\text{max}=11}^1$ of 2048 cells, the matrix $\mathcal C$ in its formulation \eqref{eq:C_1D} would have $(2^5 + \dots + 2^{11}) =4096$ lines, accounting for the levels $l=5,\dots,11$, and $4096\times 5=20480$ (unsigned) integer valued entries. Accounting for 4 bytes per integer, the memory needed to store $\mathcal C$ is $80$~KB.

\begin{figure}[t]
	\centering
	\includegraphics[width=0.55\textwidth]{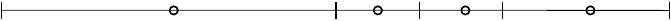}
	\caption{The 1D grid with centres $\left\{\frac{1}{4},\frac{9}{16}, \frac{11}{16}, \frac{7}{8} \right\}$ and sizes $\left\{\frac{1}{2},\frac{1}{8},\frac{1}{8},\frac{1}{4}\right\}$ is equivalently given by the set $\left\{2,12,13,7\right\}$ including the corresponding lines of the matrix $C$}\label{fig:distr_exam}
\end{figure}

				
\paragraph{Mesh representation}
Based on the above formulation, the computational grid can be represented by the (finite) sequence of line numbers of $C$ corresponding to the cells of the grid. 			
Consider, for example, the 1D mesh
\begin{equation}\label{eq:vec_FV}
	M^1=\left\{ \(\frac{1}{4},\frac 1 2\),\(\frac{9}{16},\frac 1 8\), \(\frac{11}{16},\frac 1 8\), \(\frac{7}{8}, \frac 1 4\) \right\},
\end{equation}
where each included cell is characterised by its centre $M(C)$ and its size $|D(C)|$, see also Fig.~\ref{fig:distr_exam}. The size of each cell of level $l$ is $2^{-l}$, its centre is given through the intra-level index $k$ as $\frac{2k-1}{2^{l+1}}$. The corresponding matrix $\mathcal C$ is given by~\eqref{eq_matrix_1D_example}.
The matrix lines that correspond to the mesh~\eqref{eq:vec_FV}, and in effect represent it, are given by the set
\begin{equation}\label{eq:vec_NEW}
	G_{\mathcal C}=\{2,12,13,7\}.
\end{equation}

The order in which the cells appear in $G_{\mathcal C}$ follows the \textit{Z-order} or  Lebesgue space filling curve, namely, they are order first with respect to their refinement level $l$, and then to their intra-level index $k$. This becomes most important in 2- and higher spatial dimensions.
From the implementation point of view, we use vectors of the form \eqref{eq:vec_NEW}, instead of \eqref{eq:vec_FV}, to communicate the grid between different parts of the algorithm. 
				
The following example considers the adaptation of $G_{\mathcal C}$ in the case of refinement: if one step of the refinement takes place and the cell $\(\frac{1}{4},\frac 1 2\)\in M^1$ is replaced by its children cells, the mesh cell $2\in G_{\mathcal C}$ needs to be replaced by its children cells, represented in the lines 4 and 5 of $\mathcal C$, and the grid $G_{\mathcal C}$ becomes
\[ G_{\mathcal C}^\text{new}=\{4,5,12,13,7\}. \]
				
\subsubsection*{Cell representation in 2D}
				
In dimension two, we include two additional columns in the matrix $\mathcal C$ to account for all the children, i.e. 
\begin{equation}\label{eq:C_2D}
  C=\lmatrix
  ~&   \vec{k} &~\\
  ~&  \vec{l} &~ \\
  ~&  \vec{p} &~  \\
  ~&  \vec{d_{NW}} & ~\\
  ~& \vec{d_{NE}} &~\\
  ~& \vec{d_{SW}} &~\\
  ~& \vec{d_{SE}} &~\\
	\rmatrix^T,
\end{equation}
where the intra-level index in the line vector $\vec{k}$ describes the cells of the same level in a lexicographic order first with respect to the $x$- and the $y$- coordinate of the centre of the cell. Here  $\vec{d_{NW}}$, $\vec{d_{NE}}$, $\vec{d_{SW}}$ and $\vec{d_{SE}}$ are the lines of $C^T$ corresponding to the four children cells located northwest (NW), northeast (NE), southwest (SW), and southeast (SE) with respect to the centre of the parent cell indicated by the matrix line. As in the 1D case, the line vectors $\vec{k}$, $\vec{l}$ and $\vec{p}$ include the intra-level index of the cell, its refinement level, and the line of its parent cell. Again each cell can be uniquely identified by its refinement level and its intra-level index. 

\paragraph{Size of $C$ and memory usage in 2D}
It is more instructive to describe the memory needed for storing the matrix $\mathcal C$ via an example that covers a wide class of numerical experiments. The general case follows in a similar way. 
We assume hence, a cascade of meshes $G_l^2$, with $l_\text{min}\leq l\leq l_\text{max}$ where $l_\text{min}=4$, i.e. $16\times 16$ cells are included on the coarsest mesh and $l_\text{max}=10$, i.e. $1024\times 1024$ cells are included on the finest mesh. The total number of entries of the matrix $\mathcal C$ is 9,786,112 and, with 4 bytes per input, the memory needed is approximately 39~MB. 
				
Such a memory consumption is not a strong constraint, especially since it refers to a high resolution in 2D. Confer, e.g.,~\cite{Ma-2013} where a similar mesh resolution was used to solve the first stage of a tsunami wave. Despite that, we can even more reduce the memory usage after noting that the cells on the coarsest level $l_\text{min}$ do not posses any parent cells, nor do the cells in the finest level $l_\text{max}$ posses any children cells. This means that the matrix $\mathcal C$ will need 4,194,560 less entries, and the actual memory used is reduced to approximately 22~MB. To achieve such savings in practice the columns of $\mathcal C$ are separately stored in memory blocks of different sizes.
	
An even further reduction in memory consumption can be achieved by removing the $k$- and $l$-column from the matrix $\mathcal C$. The cell level and the intra-level index can be easily recomputed on the fly as follows: if $i$ is the matrix line corresponding to a cell, its level is given by the smallest integer $l\geq l_\text{min} $ such that
$$ k = \sum_{j = l_\text{min}}^{l-1} 2^{jd} - i $$
is positive. Then, the intra-level index is also given by $k$. In the example above the memory needed reduces to approximately 11~MB this way. Note that the efficient mesh administration algorithms that we propose employ an inclusion map, see~\eqref{eq:map} below, which requires a small amount of additional memory.
		
\paragraph{Size of $C$ and memory usage in 3D}		
In the same manner we note that the memory required for a 3D case that spans from a coarse mesh of $16\times16\times16$ cells to a fine mesh of $256\times256\times256$ cells will be approximately 292~MB. These requirements can also be reduced by computing the $k$- and $l$-columns on the fly to approximately 146MB, i.e. an average of 16 or 8 bytes per mesh cell respectively. To computationally facilitate the search for neighbours, we derive an inclusion map, see  \eqref{eq:map}, which has an additional burden of 36.8 MB. These memory costs are affordable, especially if it is taken into account that the resolution $256\times256\times256$ is adequate even for challenging problems, e.g. in \cite{Dimonte-2004} such mesh resolution was deemed sufficiently accurate for the comparative study of numerical methods solving the turbulent Rayleigh-Taylor instability. Additional memory can be saved if the refinement level $l$ is stored in $\mathcal C$ as an unsigned integer of 1 byte size; this would result in a further reduction by 18 MBs.

In an even larger scale, further reduction in memory usage can be achieved by elaborate compression techniques, see e.g. \cite{Blandford2005}, but they fall beyond the scope of this paper. However, we note that the above mentioned 2- and 3D mesh resolutions are memory-wise affordable and adequate for a wide range of academic and non-academic studies.

\subsubsection*{Comparison with other grid administration techniques}

The reported memory consumption in \texttt{alugrid}~\cite{alugrid} ranges between 700-800 \texttt{bytes} per element for hexahedral elements. Due to the structure and the information stored, \texttt{alugrid} can not be directly compared to our method, but assuming a memory consumption of 146 MB, \texttt{alugrid} can store approximately the information for $2\times10^5$ cells whereas our method stores the information for $2\times10^7$ cells.

The memory consumption of \texttt{p4est} on 3D domains, as discussed in~\cite{p4estSA}, is 24 bytes per octant. Our proposed method on the other hand necessitates 16 bytes per octant. However, in our case we store the information of the full tree rather than the current mesh conformation as done in \texttt{p4est}. In a more generally setting, \texttt{p4est} aims to the scalable parallelism to thousands of processors, whereas our approach emphasizes on the simplicity, ease of implementation, and portability of the code.

		
\subsection{Siblings and neighbours}\label{sec:mat:op}
From the information stored in the matrix $\mathcal C$, the parent and the children cells of every cell are directly accessible. However, in practical situations further information is needed. In particular, the identification of the siblings and the neighbours of a given cell is important in e.g. the coarsening procedure or the computation of the numerical fluxes.				
In the following we describe these operations, based on the information stored in the matrix~$\mathcal C$. 

\subsubsection*{Siblings}				
Let a cell $C$ be represented in the line $i$ of the matrix $\mathcal C$. Then, its parent cell can be found in the line $m=\mathcal C(i,3)$, and its group of siblings in the lines $\mathcal C(m,4)$ through $\mathcal C(m,5)$ in 1D, and $\mathcal C(m,4)$ through $\mathcal C(m,7)$ in 2D.
					
To identify the relative position of a cell among its siblings, we use the ordering with respect to the (intra-level) index $k$. In 2D in particular, a cell is a \textit{south-}(S), \textit{north-}(N), \textit{west-}(W), or \textit{east-}(E) sibling according to the rules 
\begin{equation} \label{eq:siblings}
	\begin{cases}
		\text{W-cell}, & \text{if } k \text{ is odd}  \\
		\text{E-cell}, & \text{if } k \text{ is even} 
	\end{cases}, 
	\qquad \text{and} \qquad
	\begin{cases}
		\text{N-cell}, & \text{if } \left \lfloor{\frac{k}{2^{l}}}\right \rfloor \text{is odd}  \\
		\text{S-cell}, & \text{if } \left \lfloor{\frac{k}{2^{l}}}\right \rfloor \text{is even} 
	\end{cases}. 
\end{equation}

					
				
\subsubsection*{Neighbours}
To efficiently identify the neighbours of a cell, we employ the regularity of the mesh. For a RSM $G$ of smoothness $m_r$, and a cell $C\in G$ of level $L(C)=l$ and intra-level index $k$ we identify its neighbours by distinguishing the following cases:
\paragraph{Neighbours of the same level}
The neighbours  of $C$ on the uniform mesh  $G_l^d$ can be easily found using the intra-level index $k$. Neighbours to the left and right in dimension one, and to the W- and E- in dimension two, have respectively indexes $k\mp 1$, while neighbours to the S- and N- in dimension two have indexes $k\mp 2^{l-1}$. Thus they are represented in the corresponding $\mp 1$ and $\mp 2^{l-1}$ lines of the matrix $C$ relative to the cell $C$. 
The \textit{same level neighbours} however are not necessarily cells of the grid $G$. Nevertheless the identification of these $2$ cells in 1D, or $4$ cells in 2D is crucial in the neighbour finding process.
					
\paragraph{Neighbours of lower levels} 
If $N$ is a same level neighbour of $C$, but not part of the actual grid, i.e. $N\notin G$, then an ancestor $A$ of $N$ could be included in the current mesh $G$. If there is such an ancestor $A\in G$, then $A$ is also a neighbour of $C$ due to the grid structure. In order to find all possible neighbours, we check $m_r$ generations of ancestors for inclusion in the current grid. These ancestors can be identified by iteratively using the parent-cell entry $m$ of the matrix $C$ and switching to the corresponding line.
This means that for at most $2m_r$ cells in 1D and $4 m_r$ cells in 2D it needs to be examined if they are included in $G$. 
					
\paragraph{Neighbours of higher levels} 
As before, let $N$ be a same-level neighbour of the cell $C$.  If \textit{neither $N$ nor any of its ancestors} is included in $G$, we look for neighbours of $C$ among the descendants of $N$. Once again, $m_r$ generations have to be screened. 
In the algorithm we propose, we exploit the relative position of $N$, e.g. assume that $N$ is an E-neighbour of $C$ in 2D and proceed as follows:
\begin{quote}
	starting with a queue containing only $N$ we iterate through the queue by checking each entry for inclusion in $G$. If it is included, we have found a neighbour, otherwise we add the NE- and SE-child of the entry to the queue for the next iteration step. 
\end{quote}
This way all neighbours of $C$, among the descendants of $N$, can be found in at most $m_r +1$ iteration steps.
					
For an efficient computation of all the neighbours of a cell $C\in G$, we propose to compute all same level neighbours first, and afterwards 
check the same level neighbours themselves, then their ancestors, and last their descendants for being included in the grid $G$. To allow for computationally inexpensive checks for inclusion of a cell in the current mesh, we derive an inclusion map from the grid: For $G=\{C_1, \dots, C_N\}$ we store
\begin{equation}\label{eq:map}
  m_{G}(C)=\begin{cases}k & \text{if } C=C_k\in G \\ 0 & \text{if }C\notin G\end{cases}
  \end{equation}
  for every cell $C$ decoded by the matrix $\mathcal{C}$. This map is also used for the handling of approximate functions in FV schemes. Note that the additional memory requirements are minor, e.g.,  in the 3D example in Section~\ref{sec:cell:repr} approximately 36.8 MB are required.
					
%

\subsection{Refinement/Coarsening}\label{sec:ref.coars}
					
If the mesh is used for the numerical solution of PDEs, we employ \change{monitor functions that assign values to the cells}{the monitor function $M$, which assigns a non-negative value to each cell of the grid} and accordingly mark\change{ them}{s cells} for refinement and coarsening. Typically, the marking is decided by two threshold values $\theta_\text{coars}<\theta_\text{refin}$. The cells for which the monitor function is below $\theta_\text{coars}$ are marked for coarsening, whereas the cells where the value of the monitor function is greater $\theta_\text{refin}$  are marked for refinement.
					
By $G_{\text{ref}(M)}$ we denote the mesh $G$ where all cells $M \subset G$ have been refined once, i.e. replaced by their children cells, similarly $G_{\text{coars}(M)}$ denotes the mesh $G$ where all cells $M \subset G$ have been coarsened once, i.e. replaced by their corresponding parent cell.
In what follows, we assume that after evaluation of the mesh by the monitor function, we have identified the subsets $M_r$ and $M_c\subset G$ that include cells which are marked for refinement and coarsening, respectively.
					
\paragraph{Strong refinement and weak coarsening}
Note that the meshes $G_{\text{ref}(M_r)}$, and $G_{\text{coars}(M_c)}$ might not satisfy the grid regularity condition \eqref{eq:grid_reg}. Hence $M_r$ and $M_c$ have to be adjusted according to the structure of $G$. In practice, we set higher priority on refinement and lower priority on coarsening; we aim for \textit{strong refinement} and \textit{weak coarsening}. That is, we possibly refine more cells than initially marked for refinement and coarsen less cells than initially marked for coarsening. In effect, the actually refined and coarsened meshes are given by $G_{\text{ref}(M_R)}$ for $M_R\supseteq M_r$ and $G_{\text{coars}(M_C)}$ for $M_C\subseteq M_c$, respectively.
					
For each cell $C_m$ we denote the \textit{dependency sets} by
\begin{align*} D_r(C_m)&= \{ \text{smallest set }Q \subseteq G : \ C_m \in Q, \quad G_{\text{ref}(Q)} \text{ satisfies \eqref{eq:grid_reg}}\},\\
	D_c(C_m) & = \{ \text{smallest set }Q \subseteq G : \ C_m \in Q, \quad G_{\text{coars}(Q)} \text{ satisfies \eqref{eq:grid_reg}} \}. 
\end{align*}
Given these definitions we conduct strong refinement and weak coarsening using
$$M_R = \bigcup_{C_r \in M_r} D_r(C_r), \text{ and } M_C = \{C_c \in M_c \,| \	D_c(C_c) \subseteq M_c\}.$$
				 
\paragraph{Algorithm}
We conduct both, strong refinement and weak coarsening by starting with $M_R \leftarrow M_r$, $M_C \leftarrow M_c$ and iterating through $M_R$, and $M_C$ from the highest to the lowest level of the cells. If a cell $C_r \in M_R$ in the refinement process has a neighbour $N$ with $L(C_r) - L(N) = m_r$, we mark the neighbour for refinement, i.e. we add $N$ to the refinement set $M_R$. If a cell $C_c \in M_C$ in the coarsening process has a neighbour $N$ with $ L(N) - L(C_c) = m_r$, and $N \notin M_R$, we remove $C_c$ from the coarsening set $M_C$. We refer to the Algorithms~\ref{alg:refin} and~\ref{alg:coars} for further details.
				 
In numerical computations, we perform a \textit{mesh update} of a grid $G$ as follows: we first compute the monitor function, mark cells for refinement and coarsening and deduce the sets $M_r$ and $\tilde M_c$ using threshold values as described above. Then we conduct strong refinement using the effictive refinement set $M_R$. Since the refined mesh might not include all the cells that were initially marked for coarsening, i.e. possibly $\tilde M_c \nsubseteq G_{\text{ref}(M_R)}$, we afterwards conduct weak coarsening using the updated set $M_c = \tilde M_c \backslash M_R \subseteq G_{\text{ref}(M_R)}$.

\begin{algorithm}
	\caption{Strong refinement. Mesh cells are marked for refinement using an indicator function and then sorted by their refinement level in $S_{l_\text{min}},\dots,S_{l_\text{max}}$. In the following iterations sorted by refinement level, conflicts, due to condition \eqref{eq:grid_reg}, are resolved by additionally refining neighbours causing these conflicts.}\label{alg:refin}
	\begin{algorithmic}[1]
		\STATE \textbf{Input:} the current mesh $G^\text{old}$ and the numerical solution $U^\text{old}$
		\STATE Initialize the new mesh $G^\text{new} \leftarrow  G^\text{old}$
		\STATE Initialize $S_{l_\text{min}},\dots,S_{l_\text{max}}\leftarrow \emptyset$
                \FOR{$C\in G^\text{old}$}
		\IF{$M(C)>\theta_\text{refin}$ and $L(C) < l_\text{max} $ } 
		\STATE $S_{L(C)} \leftarrow S_{L(C)} \cup C$
		\ENDIF
		\ENDFOR
		
		\FOR{$l=l_\text{max}-1,\dots, l_\text{min}+m_r$}
		\FOR{$C\in S_l$}
		\STATE find the set of neighbours $N(C)$ of $C$
		\FOR{$ C_N \in N(C)$}
                \IF{$l-L(C_N)=m_r$ and $ C_N\not\in S_{l-m_r}$}
                \STATE $S_{l-m_r}\leftarrow S_{l-m_r} \cup \left\{ C_N\right\}$
		\ENDIF
		\ENDFOR
		\ENDFOR
		\ENDFOR
		
		\FOR{$C \in \bigcup_{l=l_\text{min}}^{l_\text{max}-1} S_l$}
		\STATE identify the set of children $C_c$ of $C$
		\STATE $\mathcal M^\text{new} \leftarrow \mathcal M^\text{new}\setminus \left\{ C \right\}$
		\STATE $\mathcal M^\text{new} \leftarrow \mathcal M^\text{new}\cup C_c$
		\STATE project the numerical solution from $U^\text{old}\big|_{C}$ to $U^\text{new}\big|_{C_c}$
		\ENDFOR
		\STATE \textbf{Output:} the refined mesh $G^\text{new}$ and the updated numerical solution $U^\text{new}$
	\end{algorithmic}
\end{algorithm}

\begin{algorithm}
	\caption{Weak coarsening. Mesh cells are marked for coarsening using an indicator function and then sorted by their refinement level in $S_{l_\text{min}},\dots,S_{l_\text{max}}$. In the following iterations from high to low refinement levels, coarsening of a mesh cell only takes place if all sibling cells are on the same refinement level and marked for coarsening and no conflicts due to \eqref{eq:grid_reg} are detected.}\label{alg:coars}
	\begin{algorithmic}[1]
		\STATE \textbf{Input:} the current mesh $G^\text{old}$ and the numerical solution $U^\text{old}$
		\STATE Initialize the new mesh $G^\text{new} \leftarrow G^\text{old}$
		\STATE Initialize the new numerical solution $U^\text{new} \leftarrow U^\text{old}$
		\STATE Initialize $S_{l_\text{min}},\dots,S_{l_\text{max}}\leftarrow \emptyset$
		\FOR{$C \in G^\text{old}$}
		\IF{$M(C) < \theta_\text{coar}$ and $L(C)>l_\text{min}$}
		\STATE $S_{L(C)} \leftarrow S_{L(C)} \cup \left\{C \right\}$
		\ENDIF
		\ENDFOR
		\FOR{$l=l_\text{max},\dots, l_\text{min}+1$}
		\FOR{$C\in S_l$}
		\STATE identify the set of siblings $S^b(C)$ of $C$
		\IF{$S^b(C)\not\subset S_l$}
		\STATE $S_l\leftarrow S_l\setminus \left\{S^b(C)\right\}$
		\STATE \textbf{continue} for loop in line 11
		\ENDIF
		\IF{$l_\text{max}-m_r$}
		\STATE compute the set of neighbours $N(C)$ of $C$
		\FOR{$C_N\in N(C)$}
		\IF{$L(C_N)-l=m_r$ and $C_N\not\in S_{l+m_r}$}
		\STATE $S_l \leftarrow S_l\setminus\left\{S^b(C)\right\}$
		\STATE \textbf{continue} for loop in line 11
		\ENDIF
		\ENDFOR
		\ENDIF
		\STATE identify the parent cell $C_p$ of $C$
		\STATE $G^\text{new} \leftarrow G^\text{new}\setminus \left\{S^b(C) \right\}$
		\STATE $G^\text{new} \leftarrow G^\text{new}\cup \left\{C_p\right\}$
                \STATE project the numerical solution $U^\text{old}\big|_{S^b(C)}$ to $U^\text{new}\big|_{C_p}$
		\ENDFOR
		\ENDFOR
		\STATE\textbf{Output:} the coarsened mesh $G^\text{new}$ and the updated numerical solution $U^\text{new}$
	\end{algorithmic}
\end{algorithm}


\subsection{Data handling and projection in FV}\label{sec:projections}
In practice a RSM is used to handle piecewise defined functions. Let us consider a measurable function $u: (0,1)^d \rightarrow \mathbb{R}$ and its piecewise constant representation on a RSM $G$,
$$ u_G(x) = \sum_{C \in G} U_C \chi_{D(C)}(x), $$
where
\begin{equation}\label{eq:FV} 
	U_C = |C|^{-1} \int_C u(x)\, dx
\end{equation}
and where $\chi_D$ is the characteristic function of the set $D\subset \Omega$. If we assume the representation of the grid as a vector $G=\{C_1,\,C_2, \dots, C_N\}$,
such a function can be simply stored in a corresponding vector $U_G= \{U_{C_1},\,U_{C_2}, \dots, U_{C_N}\}$.

In our implementation the information of both $G$ and $U_G$ are stored in arrays of the same size that in general allow for more entries as cells are included in the grid. Initially, we allocate arrays of a generous size, which we later extend in large increments if further memory for refinement is required. The array including the information in $G$ is not ordered and allows for holes since array entries for cells removed by coarsening are ``deleted'' and set to zero. To account for newly added cells by refinement the first free memory positions being zero are used. The cell average $U_C$ in the array for $U_G$ is stored at the same index where the corresponding cell $C$ is stored in the array for $G$. To identify a cell and associated data in memory we use the inclusion map~\eqref{eq:map}. We compute cell-wise fluxes in our numerical simulations for which this approach has been efficient. However, we note that the data needs to be reordered to efficiently use common linear algebra packages such as LAPACK~\cite{lapack99}.

\paragraph{Projection to lower levels} 

Let $G^\text{low}$ be a RSM derived by a series of subsequent coarsening operations on $G$. For each $C^l\in G^\text{low}$ we define  $\mathcal{L}_G(C_l)$ to be the set that includes either $C^l$ if $C^l\in G$ or otherwise all descendants of $C^l$ that are included in $G$. The elements in $\mathcal{L}_G(C_l)$ can be easily identified using the matrix $\mathcal{C}$. We employ the formula
\begin{equation}
	U_{C^l} = \sum_{C \in \mathcal{L}_G(C^l)} \left( 2^{L(C^l)-L(C)}\right)^d\, U_C, \label{eq:project_down}
\end{equation}
to define the projected function $u_{G^\text{low}}$. This function satisfies \eqref{eq:FV} for all $C \in G^\text{low}$, hence no accuracy is lost.
				
We apply this projection after the weak coarsening procedure to update the FV solution. Further, it can be used to compute the difference (e.g. in the discrete $L^1$ norm) of two functions on different RSMs. In this case the function on the finer grid is first projected to the coarser grid before the difference is computed.
				
\paragraph{Projection to upper levels} This projection is needed to update piecewise constant functions after strong refinement. Let us consider a RSM $G^\text{up}$, which has been derived by a series of subsequent refinement operations on $G$. For each cell $C^u \in G^\text{up}$ we define by $\mathcal{U}(C^u)$ either $C^u$ itself if $C^u \in G$ or otherwise the ancestor of $C^u$ which is included in $G$. Once again, $\mathcal{U}(C^u)$ can be easily identified using $\mathcal{C}$. A simple projection for the definition of $u_{G^\text{up}}$ would be the choice
\begin{equation}
	U_{C^u} = U_{\mathcal{U}(C^u)}. \label{eq:project_up}
\end{equation}
Higher order projections, which are not considered in this work, make use of reconstructions that take into account also the neighbour cells of $\mathcal{U}(C^u)$, these include minmod, CWENO, and more~\cite{harten1987uniformly, shu1988efficient, sempliceAdaptiveMeshRefinement2016}. 

\section{Numerical experiments}\label{sec:num:exper}
We present three numerical applications of the mesh administration technique and the AMR method. The first one is a generic experiment without any physical or biological interpretation.  In this, we demonstrate the properties of the mesh administration and the AMR techniques using predefined monitor functions to drive the grid reconstruction. The second one is a 2D explosion experiment using the 2D Euler system showing that our methods can efficiently capture the more important regions of the phenomenon. The third one is a biological application; a 2D cancer invasion model exhibiting highly dynamic and complex solutions.
The computer programs employed are implemented in Fortran 2008 (mesh administration and numerical schemes) and Python 3 (visualization and error computation) and are available upon request.
\subsection{Generic experiment}\label{sec:generic}
In a first test, we consider a 2D RSM with refinement levels set to $l_\text{min}=5$ and $l_\text{max}=7$ and mesh regularity $m_r=1$. We start with a uniform mesh on the lowest level, i.e. $G = G_5^2$, and consider the time dependent monitor function
\begin{equation}\label{eq:M1}
	M_1(\mathbf x,t) = e^{-100\,\left\|\mathbf x-(0.1+t) (1,~1)^T\right\|_2^2}. 
\end{equation}
Using the small time step $\Delta t = 0.005$ and starting at $t=0$, we advance in time and perform strong refinement (without coarsening) once at every time step with threshold $\theta_\text{ref}= 0.8$ until the time $t = 0.8$ is reached.
	
\begin{figure}[t] 
    \begin{tikzpicture}
	\begin{groupplot}[
          /tikz/mark size=1.5pt,
          group style={group name=my plots, group size=4 by 1, horizontal sep=.8cm},
          xmin=0, xmax=1, ymin=0, ymax=1,
          ticklabel style = {font=\scriptsize},
          axis background/.style={fill=white},
          width=.3\linewidth, height=.3\linewidth, 
          ]
          \nextgroupplot[title = {$t=0.4$}]
          \addplot [forget plot] graphics [xmin=0, xmax=1, ymin=0, ymax=1] {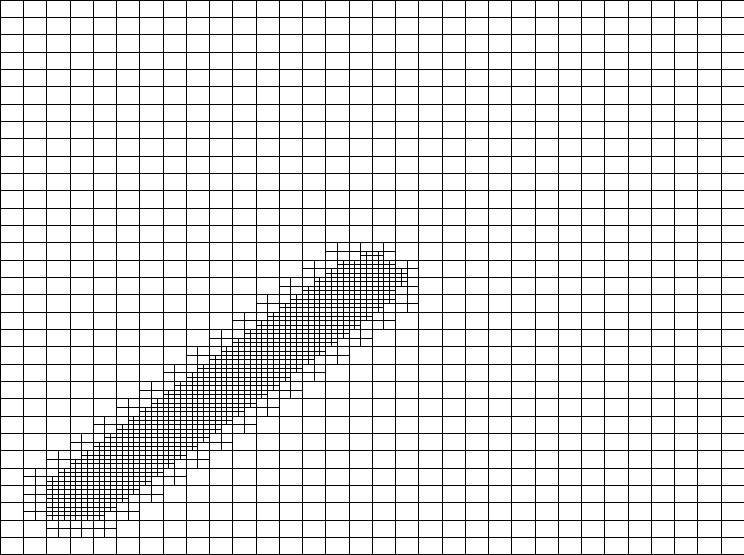};
          \nextgroupplot[title = {$t=0.8$}]
          \addplot [forget plot] graphics [xmin=0, xmax=1, ymin=0, ymax=1] {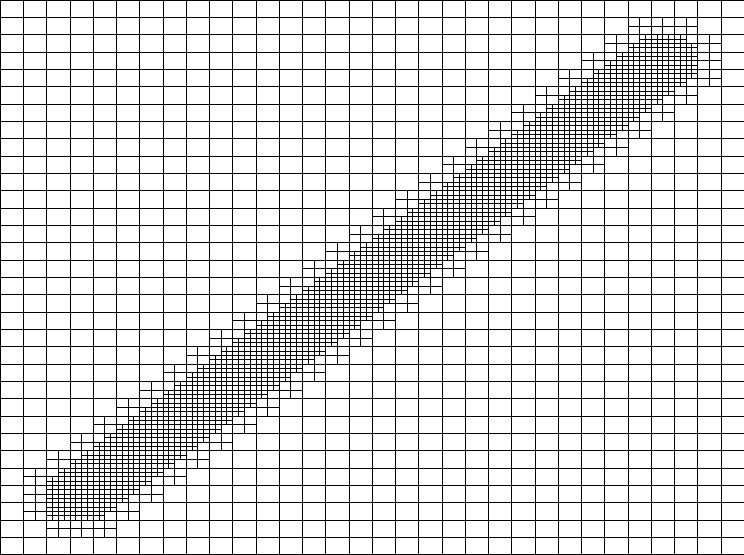};
          \nextgroupplot[title = {1 $\times$ coarsening}]
          \addplot [forget plot] graphics [xmin=0, xmax=1, ymin=0, ymax=1] {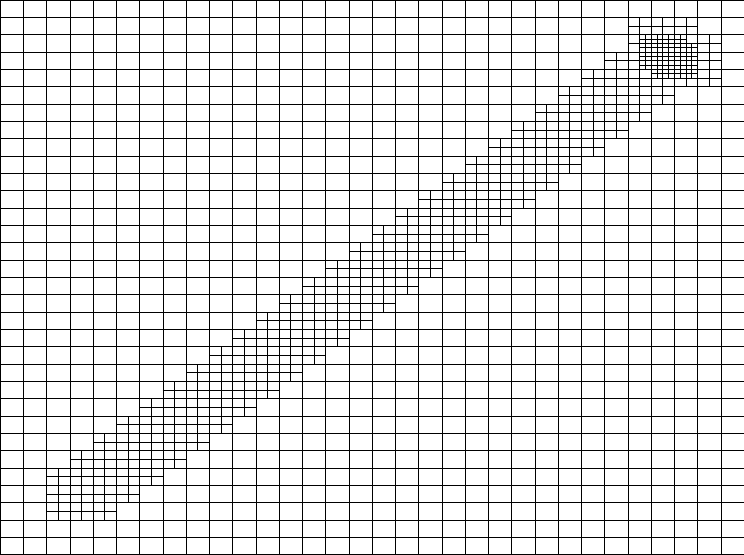};
          \nextgroupplot[title = {2 $\times$ coarsening}]
          \addplot [forget plot] graphics [xmin=0, xmax=1, ymin=0, ymax=1] {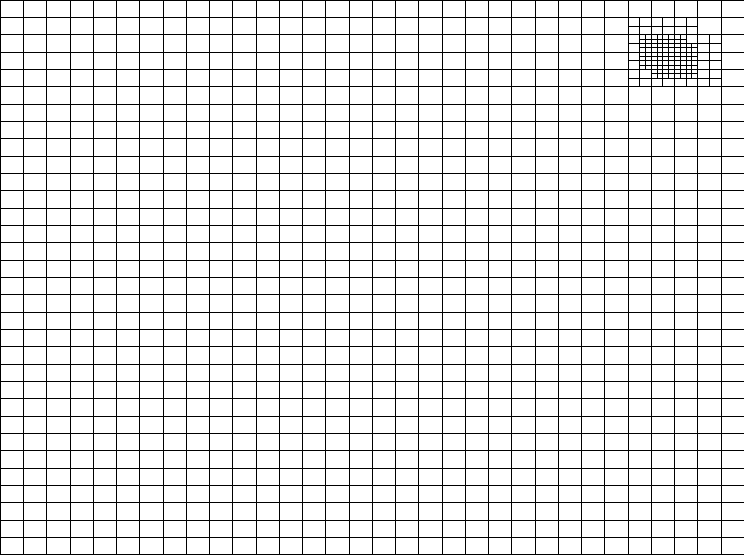};
        \end{groupplot}
      \end{tikzpicture}%
\caption{Diagonal movement of the monitor function~\eqref{eq:M1} with refinement at intermediate and final time and after one and two subsequent coarsening steps.}\label{fig:moveAndCoarsen}
\end{figure}
				
The resulting meshes at $t=0.4$ and at the final time are shown in Figure~\ref{fig:moveAndCoarsen}. The movement of the Gaussian $M_1(\mathbf x,t)$ has left a trace on the grid. On the trace of the Gaussian, the cells have been refined to the maximal level $7$. Cells of level $6$ are only visible on the transition between the finer and the coarser grid and are a result of the mesh regularity.
We then coarsen the mesh twice using $\theta_\text{coars}= 0.8$ without evolving the monitor function further in time. Figure \ref{fig:moveAndCoarsen} exhibits the mesh after performing each weak coarsening procedure. The trace of the movement is completely coarsened after the second coarsening.
				
In a second test we employ another 2D RSM with $l_\text{min}=3$ and $l_\text{max}=7$ and mesh regularity $m_r=1$. We start again with a uniform mesh on the lowest level, i.e. $G = G_3^2$, and consider the time dependent monitor functions
\begin{align}
	M_2(\mathbf x, t) &=  e^{-100\,\left\|\mathbf x- 0.9 \(\cos( 0.5\, \pi t ),~\sin(0.5 \, \pi t)\)^T\right\|_2^2}, \label{eq:M2} \\
	M_3(\mathbf x, t) &= e^{-100\,\left \|\mathbf x- 0.9 \(\cos(\pi (1 - 0.5 \, t)),~\sin( \pi (1 - 0.5 \, t))\)^T\right\|_2^2} + M_2(t), \label{eq:M3} \\
	M_4(\mathbf x,t) & = \begin{cases}1, & 0.07  < \|x\|_2 - 0.5 \, t< 0.1 \text{ or  } 0.98 <  x_1 + 1.5\, t < 1.02,\\ 0, &\text{otherwise} \label{eq:M4}
	\end{cases}.
\end{align}
Starting at $t=0$ and using the time step $\Delta t = 0.005$, we evolve in time and perform strong refinement and weak coarsening once at every time step. 

\begin{figure}[t] 
    \begin{tikzpicture}
	\begin{groupplot}[
          /tikz/mark size=1.5pt,
          group style={group name=my plots, group size=4 by 1, horizontal sep=.8cm},
          xmin=0, xmax=1, ymin=0, ymax=1,
          ticklabel style = {font=\scriptsize},
          axis background/.style={fill=white},
          width=.3\linewidth, height=.3\linewidth, 
          ]
          \nextgroupplot[title = {$t=0.05$}]
          \addplot [forget plot] graphics [xmin=0, xmax=1, ymin=0, ymax=1] {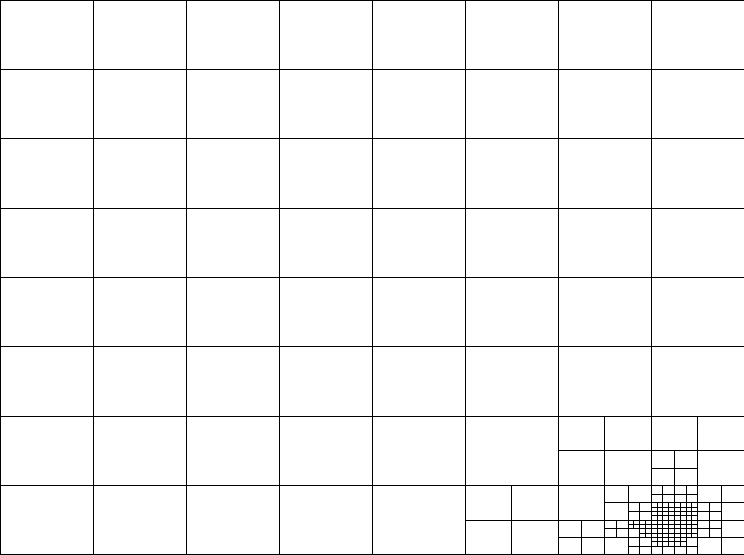};
          \nextgroupplot[title = {$t=0.25$}]
          \addplot [forget plot] graphics [xmin=0, xmax=1, ymin=0, ymax=1] {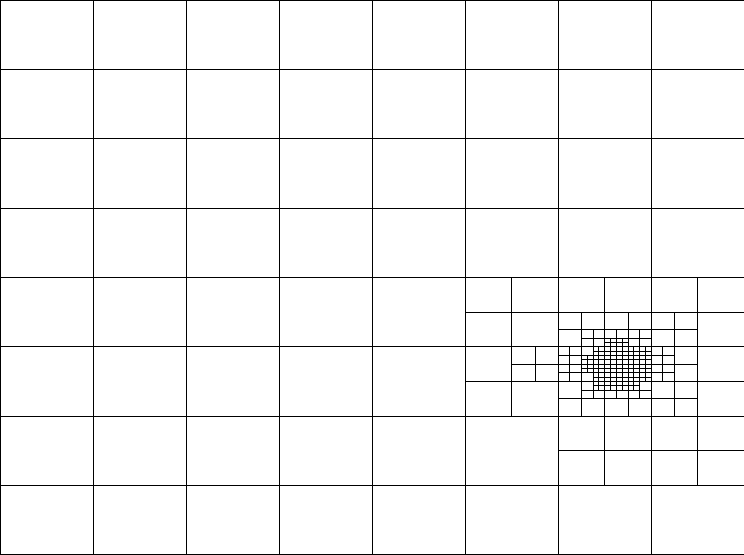};
          \nextgroupplot[title = {$t=0.5$}]
          \addplot [forget plot] graphics [xmin=0, xmax=1, ymin=0, ymax=1] {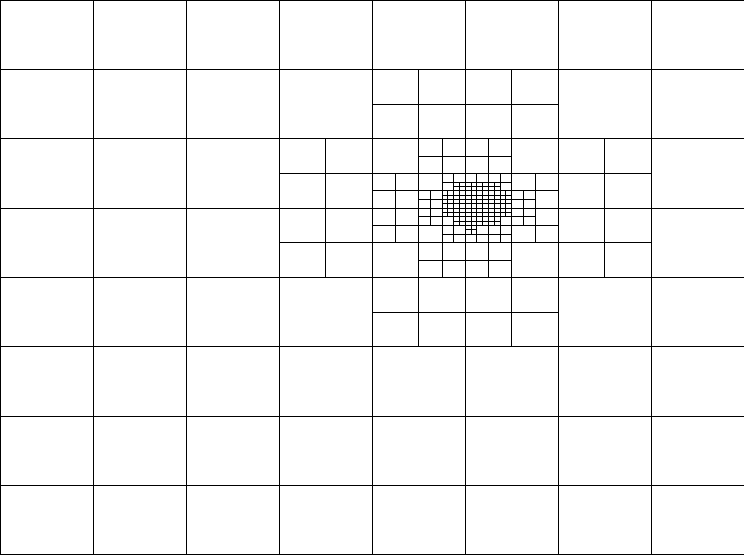};
          \nextgroupplot[title = {$t=0.8$}]
          \addplot [forget plot] graphics [xmin=0, xmax=1, ymin=0, ymax=1] {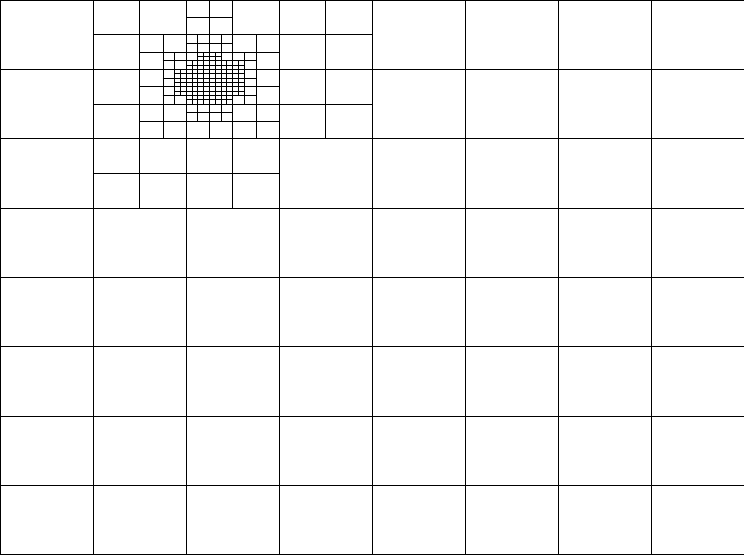};
        \end{groupplot}
      \end{tikzpicture}%
\caption{Circular movement of the monitor function \eqref{eq:M2}. Refinement and coarsening take place with the time evolution.}\label{fig:movingGaussianNoTail}
\end{figure}

In case of the monitor function $M_2$, and $\theta_\text{ref}=\theta_\text{coars}= 0.8$, we see in Figure~\ref{fig:movingGaussianNoTail} that the circular movement of the Gaussian is not memorized by the mesh. We can observe a symmetric stepwise decrease of the cell levels when moving away from the centre of the Gaussian, which is refined to the maximal level.  
					
\begin{figure}[t] 
    \begin{tikzpicture}
	\begin{groupplot}[
          /tikz/mark size=1.5pt,
          group style={group name=my plots, group size=4 by 1, horizontal sep=.8cm},
          xmin=0, xmax=1, ymin=0, ymax=1,
          ticklabel style = {font=\scriptsize},
          axis background/.style={fill=white},
          width=.3\linewidth, height=.3\linewidth, 
          ]
          \nextgroupplot[title = {$t=0.1$}]
          \addplot [forget plot] graphics [xmin=0, xmax=1, ymin=0, ymax=1] {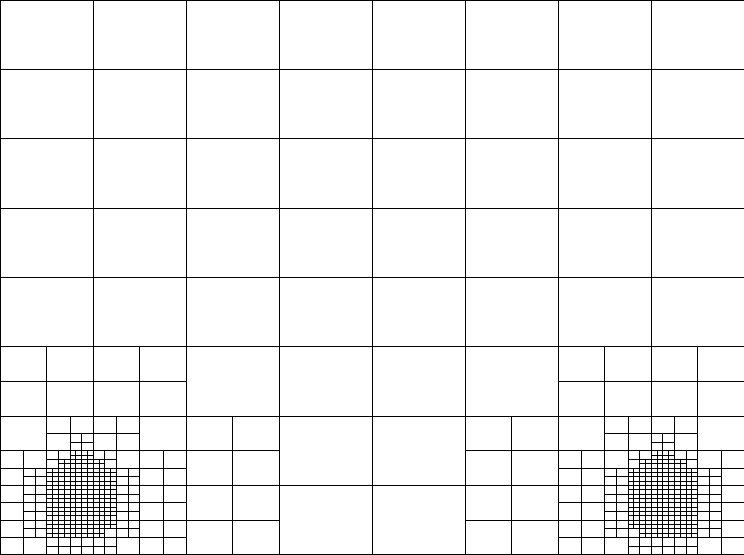};
          \nextgroupplot[title = {$t=0.25$}]
          \addplot [forget plot] graphics [xmin=0, xmax=1, ymin=0, ymax=1] {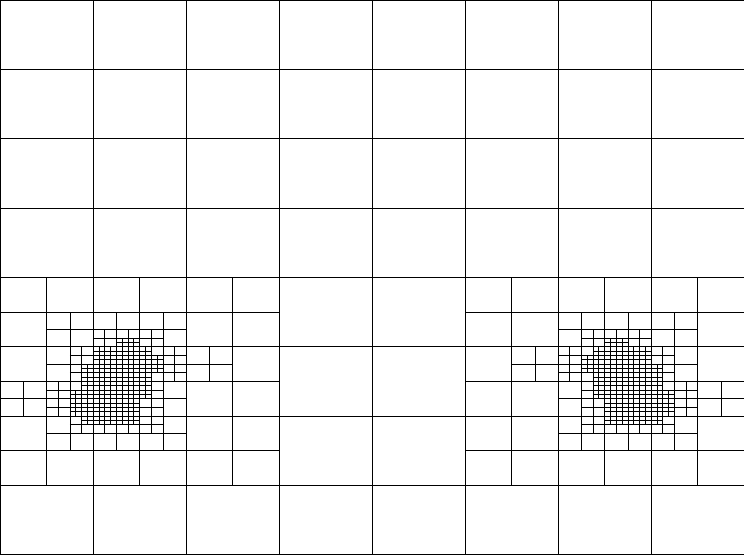};
          \nextgroupplot[title = {$t=0.5$}]
          \addplot [forget plot] graphics [xmin=0, xmax=1, ymin=0, ymax=1]  {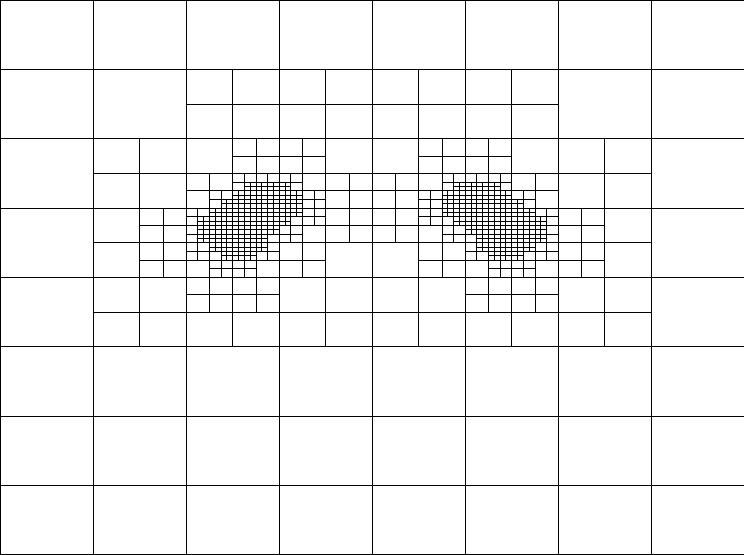};
          \nextgroupplot[title = {$t=0.65$}]
          \addplot [forget plot] graphics [xmin=0, xmax=1, ymin=0, ymax=1]  {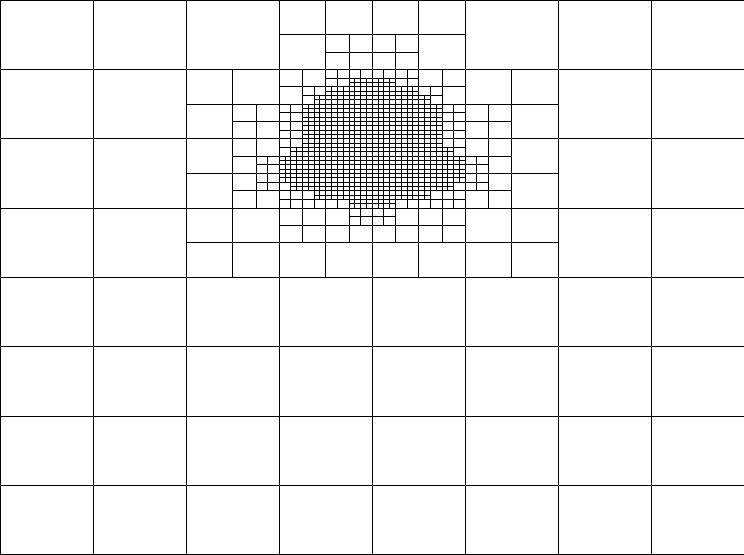};
        \end{groupplot}
      \end{tikzpicture}%
      \caption{Diagonally moving monitor functions meet. Note their ``tails'' obtained by setting the refinement and coarsening parameters $\theta_\text{refin}$, and $\theta_\text{coars}$ to different values.}\label{fig:movingGaussianTwo}
\end{figure}

Figure~\ref{fig:movingGaussianTwo} shows the same experiment using the monitor $M_3$ and the distinct thresholds $\theta_\text{ref}=0.8, \theta_\text{coars}= 0.3$. A particular memory effect of the mesh structure can be observed in the form of ``tails'' following two traveling Gaussians. The appearance of the tails in this experiment is due to the different refinement and coarsening thresholds and the particular monitor function (cf. Figure~\ref{fig:movingGaussianNoTail}). Depending on the problem under consideration, e.g. in a combined shock and rarefaction wave resolved by a FV, such property might be beneficial for the numerical investigations. 
						
\begin{figure}[t] 
    \begin{tikzpicture}
	\begin{groupplot}[
          /tikz/mark size=1.5pt,
          group style={group name=my plots, group size=4 by 1, horizontal sep=.8cm},
          xmin=0, xmax=1, ymin=0, ymax=1,
          ticklabel style = {font=\scriptsize},
          axis background/.style={fill=white},
          width=.3\linewidth, height=.3\linewidth, 
          ]
          \nextgroupplot[title = {$t=0.05$}]
          \addplot [forget plot] graphics [xmin=0, xmax=1, ymin=0, ymax=1] {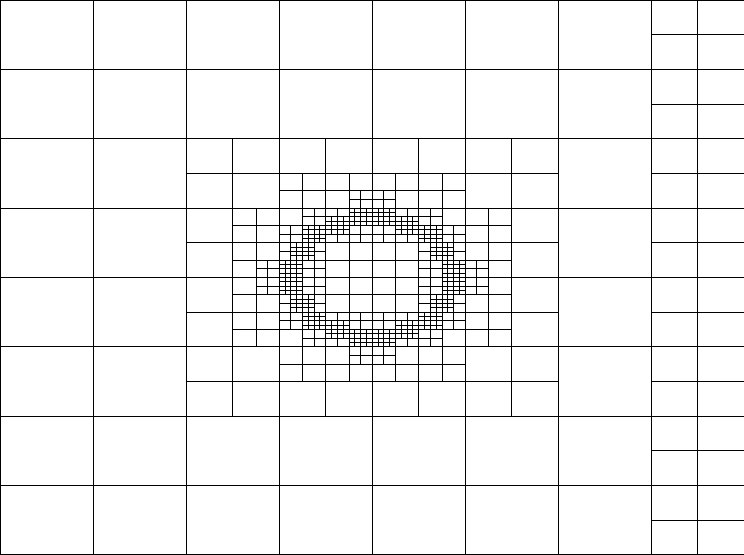};
          \nextgroupplot[title = {$t=0.18$}]
          \addplot [forget plot] graphics [xmin=0, xmax=1, ymin=0, ymax=1] {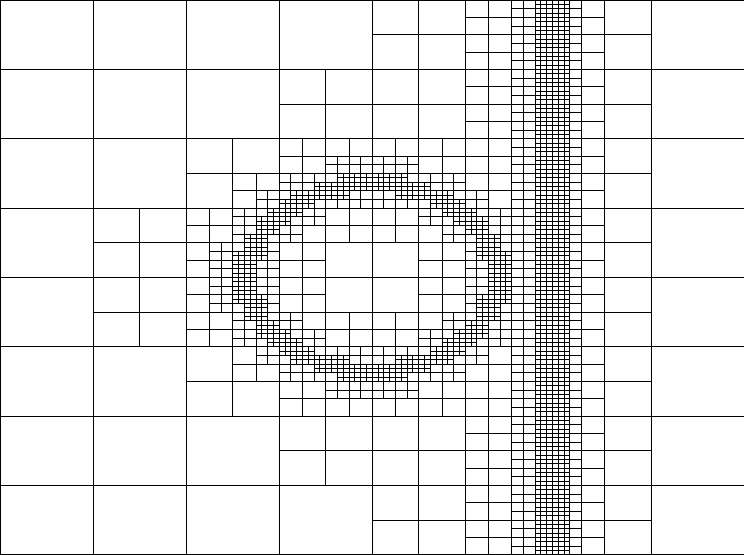};
          \nextgroupplot[title = {$t=0.38$}]
          \addplot [forget plot] graphics [xmin=0, xmax=1, ymin=0, ymax=1]  {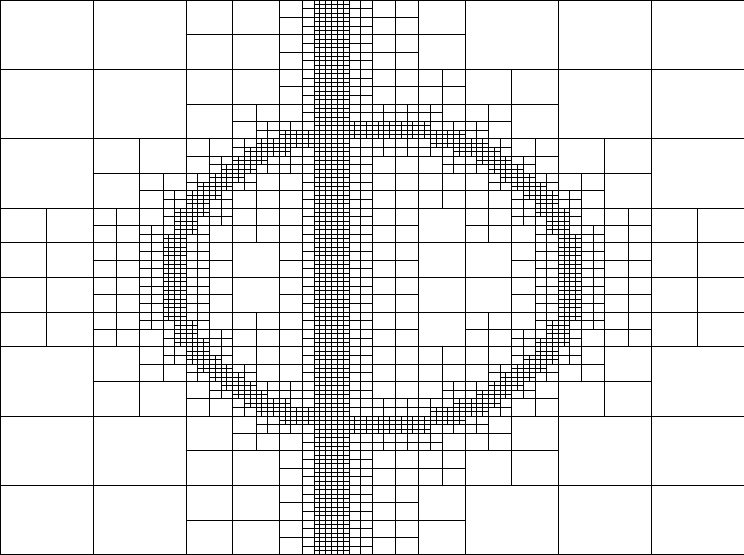};
          \nextgroupplot[title = {$t=0.55$}]
          \addplot [forget plot] graphics [xmin=0, xmax=1, ymin=0, ymax=1]  {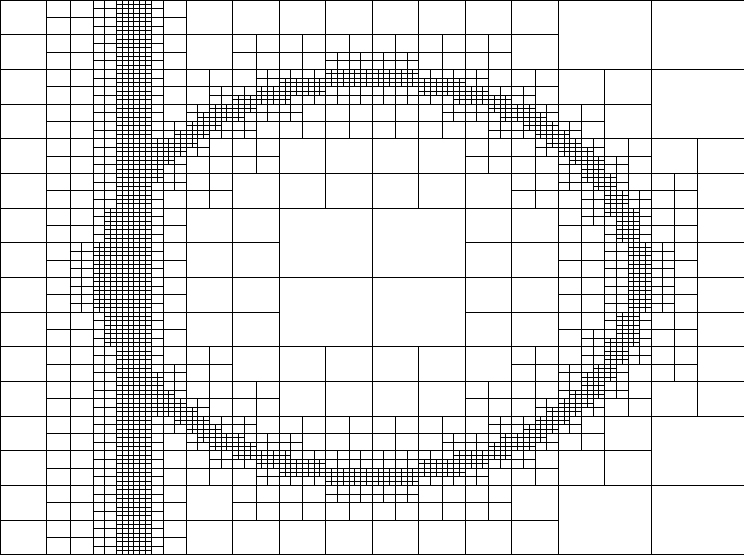};
        \end{groupplot}
      \end{tikzpicture}%
      \caption{Ring of increasing diameter crossed by a plane wave used as monitor function.}\label{fig:ring}
\end{figure}

In Figure~\ref{fig:ring} we chose as monitor function a ring of increasing diameter intersected by a plane wave given in \eqref{eq:M4} and $\theta_\text{ref}=\theta_\text{coars}= 0.5$. As the horizontally moving wave intersects the ring, cells in the tail are marked for coarsening but due to the strategy we employ and the presence of the ring some of these cells are not coarsened.  
						
\subsection{Euler equations}\label{sec:euler}
In this section we apply our mesh administration and refinement technique in a FV method to solve the Euler System
\begin{align}
  \vec U_t + \nabla \cdot \vec F(\vec U) & = 0\quad\text{in }\Omega, \quad \vec U = (\rho, \mathbf u, E)^T, \quad \mathbf u = (u_1,u_2)^T, \\
	\vec U(\cdot, 0)             & = \vec U_0(\cdot, t)\quad\text{in }\Omega                                                      
\end{align}
on the square $\Omega=(0,1)^2$, where $\rho$, $\mathbf u$, and $E$ represent the \textit{density}, \textit{velocity}, and \textit{energy} of the fluid per unit volume. The system is equipped with transmissive boundary conditions~\cite{giles1988non}, as well as with flux functions and pressure 
\begin{equation}
  \vec F( \vec U) = \begin{pmatrix}
	\rho \mathbf u\\ \rho \mathbf u \otimes \mathbf u + p I\\ (\rho E + p) \mathbf u
	\end{pmatrix}, \quad \frac{p}{\gamma - 1} = E - \frac{\rho}{2}(u^2 + v ^2)
\end{equation}
with an assumed \textit{specific heat} ratio of $\gamma = 1.4$. 
					
To denote the FV scheme we enumerate the cells of the RSM by $G^n = \{ C_i^n,~ i=1,\dots, N^n\}$. The mesh $G^n$ as well as the number of included cells $N^n$ changes throughout the computation due to the mesh adaptation that we employ. Further, we denote the set of neighbours of a cell $C$ in $G$ by $N(C)$. Given a time grid $t^0 = 0, ~ t^n< t^{n+1} = t^{n} + \Delta t^n,~n= 0, 1, 2 \dots$, we consider the cell averages
$$ \vec{\tilde U}_i^n = \frac{1}{|C_i^n|}\int_{C_i^n} \vec U(\mathbf x,t^n)\,dx,\quad i=1,\dots,N, $$
for which the exact evolution reads
\begin{equation}\label{eq:ev_nonUniform}
  \vec{ \tilde U}_i^{n+1} =  \vec{ \tilde U}_i^n - |C_i^n|^{-1} \sum_{C_j^n \in N(C_i^n)}\int_{t^n}^{t^{n+1}} 
	\int_{\partial C_{ij}^n} \vec F (\vec U(\mathbf x, t)) \vec n_{ij}^n \,d\vec x\, dt
\end{equation}
for $i=1,\dots,N^n$. Here, $\partial C_{ij}^n$ denotes the edge between the cells $C_i^n$ and $C_j^n$ and $\vec{n}_{i,j}^n$ is the outer normal vector of $C_i^n$ pointing towards $C_j^n$. By abuse of notation we assume in \eqref{eq:ev_nonUniform} that $\tilde U_i^n$ already refers to the average after projection to $C_i^n$ and hence $\tilde U_i^n$ and $\tilde U_i^{n+1}$ are averages over the same cell. The same is assumed for approximations in the following description. Employing approximate averages, $\vec U_i^n\approx \vec{ \tilde  U}_i^n$, and a numerical flux function $\vec H$, \eqref{eq:ev_nonUniform} transitions into the FV scheme
\begin{equation}\label{eq:CLScheme}
  \vec 	U_i^{n+1} = \vec U_i^n - \Delta t^n \sum_{C_j^n \in N(C_i^n)}
	\frac{|\, \partial C_{ij}^n\, |}{|C_i^n|} \,  \vec H( \vec U_i^n, \vec U_j^n, \mathbf n_{i,j}^n),
\end{equation}
for $i=1,\dots,N^n$. A time update of the FV scheme requires the identification of neighbours as well as an evaluation of the numerical flux function for each pair of neighbours. The RSM that we use allows us to compute the relation between neighbouring cell sizes and interfaces by
$$ \frac{|\,\partial C_{ij}^n \, |}{|C_i|} = 2^{2\, L(C_i^n) - \max\{ L(C_i^n), \, L(C_j^n)\}}.$$
We use a common vector-splitting approach, see~\cite{steger1981flux}, to compute both a numerical flux function\change{}{~$\vec H( \vec U_i^n, \vec U_j^n, \vec n_{ij}^n) $ as an approximation to the mean value of the actual flux $\mathbf F$ through $C_{ij}^n$, i.e. of}
$$ \frac{1}{ \Delta t^n \,|\partial C_{ij}^n|} \int_{t^n}^{t^{n+1}} 
\int_{\partial C_{ij}^n} \vec{F}(\vec U(\mathbf x, t)) \mathbf n \,d \vec x\, dt $$
and propagation velocities $a_{ij}^n$ at each cell interface $\partial C_{ij}^n$. The latter allow us to update the time increment of our explicit scheme by the local CFL condition
\[
\Delta t^n \leq \text{CFL} \min_{1\leq i,j\leq N^n} \frac{|\partial C_{ij}^n|}{|a_{ij}^n |}.
\]
A CFL number of $0.5$ is used in our simulations.

We consider the explosion experiment (cf.~\cite{Toro.2009}) with initial condition
$$ U^0(\mathbf x) = \begin{cases}
(1, 0, 2.5)^T, & \text{if }\mathbf x \in K =\{ \mathbf x \in \mathbb{R}^2, ~ \| \mathbf x - (0.5, 0.5)^T\|_2 < 0.12\},\\
(0.125, 0, 0.25)^T, & \text{otherwise}
\end{cases}.$$ 
Using a RSM with $l_\text{min} = 7, ~l_\text{max}=9$ and starting on the coarsest grid, i.e. $G^0 = G_7^2$, we perform strong refinement and subsequent weak coarsening once after each time step of scheme~\eqref{eq:CLScheme}. Therefore we use the density gradient monitor (cf. also \cite{feistauer2003mathematical}) defined on each cell by
$$g_i^n = \max_{C_j^n \in N(C_i^n)} \frac{|\rho_j^n - \rho_i^n|}{\|M(C_i^n)-M(C_j^n)\|_2}, \quad M_i^n = \frac{g_i^n}{\max_{\{ 0\leq i\leq N^n\}} g_i^n}, $$
and the threshold values $\theta_\text{ref}=\theta_\text{coars}= 0.4$. We apply the projections~\eqref{eq:project_down} and~\eqref{eq:project_up} to update the numerical solutions.
				
\begin{figure}
   \pgfplotsset{colormap/YlOrRd-9}
	\centering 
	\begin{tabular}{ccc}
        \begin{tikzpicture}
          \begin{axis}[
          title = {$t=0$},
          /tikz/mark size=1.5pt,
          xmin=0.25, xmax=0.5, ymin=0.25, ymax=0.5,
          ticklabel style = {font=\scriptsize},
          axis background/.style={fill=white},
          width=.35\linewidth, height=.35\linewidth, 
          ]
          \addplot [forget plot] graphics [xmin=0, xmax=1, ymin=0, ymax=1] {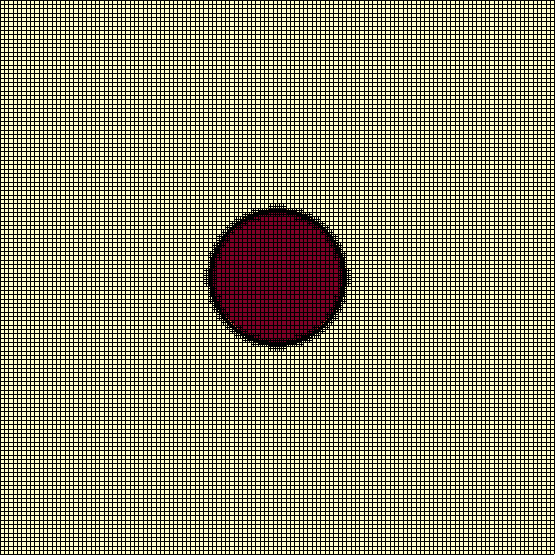};
          \end{axis}
        \end{tikzpicture}  
          &
        \begin{tikzpicture}
          \begin{axis}[
          title = {$t=0.05$},
          /tikz/mark size=1.5pt,
          xmin=0.25, xmax=0.5, ymin=0.25, ymax=0.5,
          ticklabel style = {font=\scriptsize},
          axis background/.style={fill=white},
          width=.35\linewidth, height=.35\linewidth, 
          ]
          \addplot [forget plot] graphics [xmin=0, xmax=1, ymin=0, ymax=1] {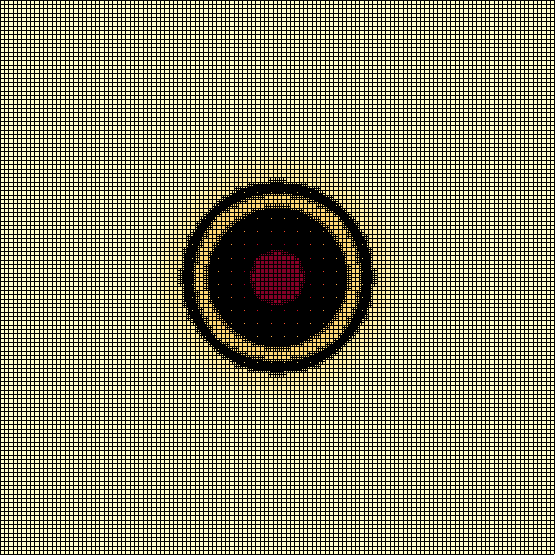};
          \end{axis}
        \end{tikzpicture} 
          &
        \begin{tikzpicture}
          \begin{axis}[
          title = {$t=0.12$},
          /tikz/mark size=1.5pt,
          xmin=0.25, xmax=0.5, ymin=0.25, ymax=0.5,
          ticklabel style = {font=\scriptsize},
          axis background/.style={fill=white},
          width=.35\linewidth, height=.35\linewidth, 
          ]
          \addplot [forget plot] graphics [xmin=0, xmax=1, ymin=0, ymax=1] {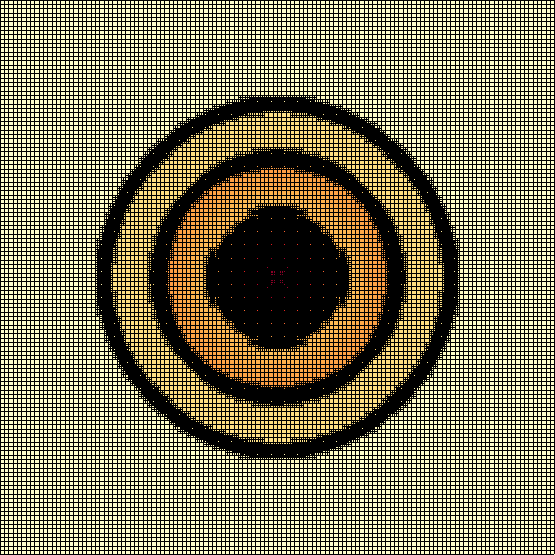};
          \end{axis}
        \end{tikzpicture} 
		\\
    \begin{tikzpicture}
	\begin{groupplot}[
          /tikz/mark size=1.5pt,
          group style={group name=my plots, group size=2 by 1, horizontal sep=1cm},
          xmin=0, xmax=1, ymin=0, ymax=1, ticklabel style = {font=\scriptsize},
          axis background/.style={fill=white}, xtick={0,1}, ytick={0,1},
          width=.2\linewidth, height=.2\linewidth,
          colorbar/width=1.5mm, colorbar, colorbar horizontal,
          colorbar style={at={(0,-0.4)},anchor=north west}
          ]
          \nextgroupplot[title={$\rho$}, point meta min=0.1, point meta max=1]
          \addplot [forget plot] graphics [xmin=0, xmax=1, ymin=0, ymax=1] {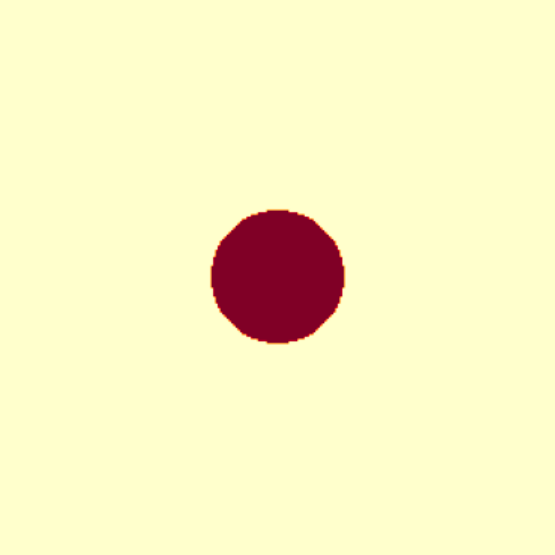};
           \nextgroupplot[title={$E$}, point meta min=0.3, point meta max=2.5]
          \addplot [forget plot] graphics [xmin=0, xmax=1, ymin=0, ymax=1] {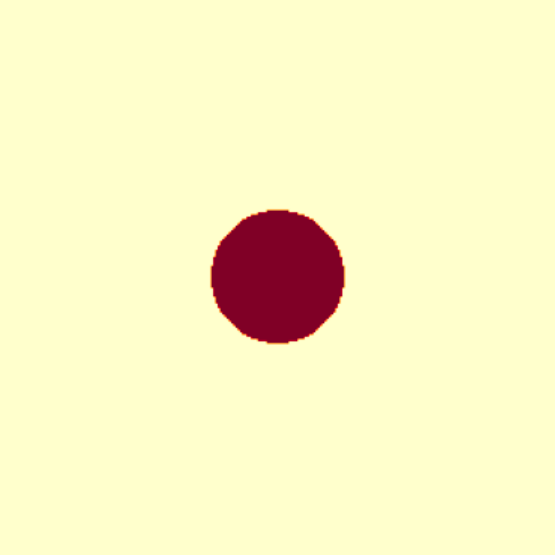};
          \end{groupplot}
      \end{tikzpicture}%
          &
            \begin{tikzpicture}
	\begin{groupplot}[
          /tikz/mark size=1.5pt,
          group style={group name=my plots, group size=2 by 1, horizontal sep=1cm},
          xmin=0, xmax=1, ymin=0, ymax=1, ticklabel style = {font=\scriptsize},
          axis background/.style={fill=white}, xtick={0,1}, ytick={0,1},
          width=.2\linewidth, height=.2\linewidth,
          colorbar/width=1.5mm, colorbar, colorbar horizontal,
          colorbar style={at={(0,-0.4)},anchor=north west}
          ]
          \nextgroupplot[title={$\rho$}, point meta min=0.1, point meta max=1]
          \addplot [forget plot] graphics [xmin=0, xmax=1, ymin=0, ymax=1] {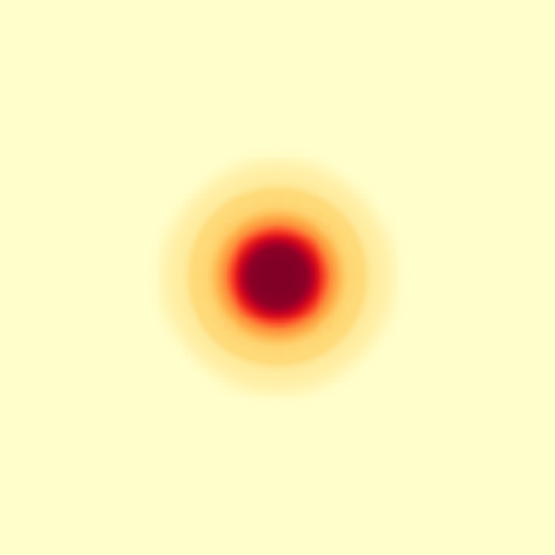};
           \nextgroupplot[title={$E$}, point meta min=0.25, point meta max=2.5]
          \addplot [forget plot] graphics [xmin=0, xmax=1, ymin=0, ymax=1] {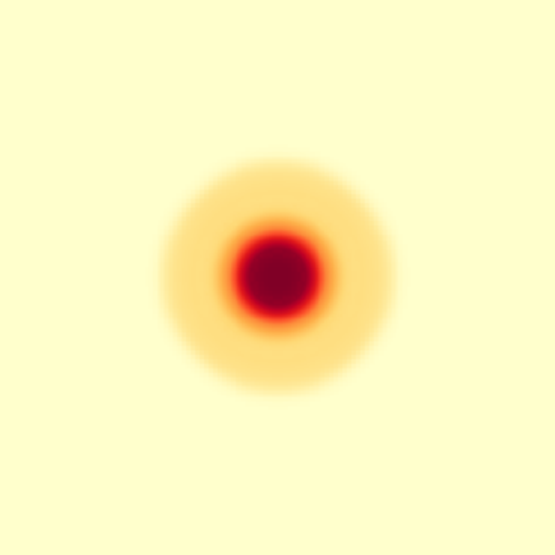};
          \end{groupplot}
      \end{tikzpicture}%
          & 
\begin{tikzpicture}
	\begin{groupplot}[
          /tikz/mark size=1.5pt,
          group style={group name=my plots, group size=2 by 1, horizontal sep=1cm},
          xmin=0, xmax=1, ymin=0, ymax=1, ticklabel style = {font=\scriptsize},
          axis background/.style={fill=white}, xtick={0,1}, ytick={0,1},
          width=.2\linewidth, height=.2\linewidth,
          colorbar/width=1.5mm, colorbar, colorbar horizontal,
          colorbar style={at={(0,-0.4)},anchor=north west}
          ]
          \nextgroupplot[title={$\rho$}, point meta min=0.12, point meta max=0.51]
          \addplot [forget plot] graphics [xmin=0, xmax=1, ymin=0, ymax=1] {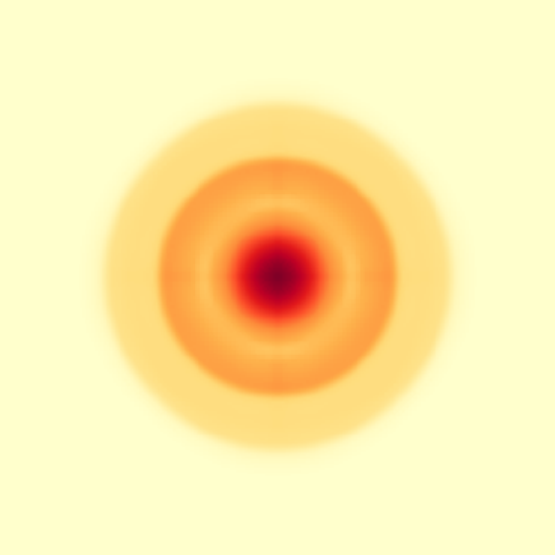};
          \nextgroupplot[title={$E$}, point meta min=0.2, point meta max=0.91,
          colorbar style = {xtick={0.3, 0.8}}]
          \addplot [forget plot] graphics [xmin=0, xmax=1, ymin=0, ymax=1] {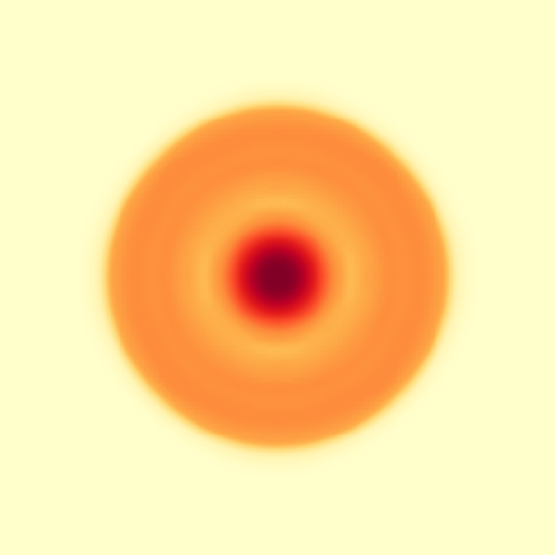};
          \end{groupplot}
      \end{tikzpicture}%
	\end{tabular}
	\caption{Euler explosion in 2D with adaptive mesh refinement. Numerical solutions at the time instances $t=0$, $t=0.05$ and $t=0.12$. Top: density and mesh in $[0.25, 0.5]^2$. Bottom: density and energy on the full domain. A shock, a contact discontinuity, and a rarefaction wave are formed and resolved on the non-uniform grid. The numerical solution employs 16,144 coarse cells ($l=7$) and 612 ($l=8$) + 1392 ($l=9$) refined cells at time $t=0$, 15,136 coarse cells ($l=7$) and 1208 ($l=8$) +~15,136 ($l=9$) refined cells at time $t=0.05$ and 13,504 coarse cells ($l=7$) and 2864 ($l=8$) +~34,624 ($l=9$) refined cells at time $t=0.12$.}\label{fig:Euler}
\end{figure}

In Figure~\ref{fig:Euler} we show the results of the simulation. We observe the formation of a circular shock wave travelling in an outward direction, followed by a contact discontinuity that also moves outwards, and by a rarefaction wave that moves inwards. We see that all three waves are properly resolved by the mesh, most notably the contact discontinuity, despite the lower order of accuracy of the numerical scheme that we have adopted. Note also that the transition areas between the waves are resolved by the highest refinement level.
	
\subsection{Tumor growth model}\label{sec:cancer}
				
\textit{Cancer metastasis} starts with the invasion of the local \textit{extracellular matrix} (ECM) by the cancer cells of the primary tumour. Then, it evolves as follows: the cancer cells proliferate and migrate until they find neighbouring blood vessels. They then \textit{intravasate} (enter the blood stream) and travel with the blood flow. At a secondary location of the organism they \textit{extravasate} (exit from the blood stream) to a new organ and engender new tumours. For that reason, the invasion of the ECM is considered to be one of the ``\textit{hallmarks of cancer}''~\cite{hanahan2011hallmarks}.
				
The model we address in this work is a 2D, macroscopic, deterministic, ARD system that describes the densities of the cancer cells $c$ as the primer unknown variable, the density $v$ of the \textit{collagen} on which cancer cells adhere and move as the main component of the \textit{extracellular matrix} (ECM), and the density of a generic enzyme $m$ of the \textit{matrix metalloproteinases} (MMPs) family that is secreted by the cancer cells and is responsible for the degradation of the ECM.

There is a wide variety of  cancer invasion models of this type that has been proposed in the literature e.g. \cite{Chaplain.2005, Maini2008, Sfak-Kol-Hel-Luk.2016,  Hel-Kol-Sfak.2015, Preziosi2015, Preziosi.2003,Stinner.2015, Czochra.2012,Johnston.2010};  we employ here a model similar to \cite{Anderson.2000} due primarily to its simplicity. This model reads
\begin{align} \label{eq:cancer:1}
	\vec U_t + \nabla \cdot \vec F( \vec U, \nabla \vec U) = \vec S( \vec U), \quad \vec U = (c, v, m)^T \quad                             & \text{in }\Omega,          \\
	\vec U(\cdot, 0) = \vec U_0(\cdot, t)\quad\text{in }\Omega, \quad\frac{\partial \vec U}{\partial \bf n} = 0\quad & \text{on }\partial \Omega, 
\end{align}
where we have assumed homogeneous Neumann boundary conditions on the square $\Omega=(0,1)^2$ and flux function and source term given by
\begin{equation}\label{eq:cancer:2}
	\vec F(\vec U, \nabla \vec U) = \begin{pmatrix}
	\chi c \nabla v - D_c \nabla c\\ 0 \\
	-D_m \nabla m 
	\end{pmatrix}, \quad 
	\vec S( \vec U) = \begin{pmatrix}
	\mu c(1-c) \\
	-\delta v m \\
	\alpha c - \beta m
	\end{pmatrix}.
\end{equation}
				
One of the reasons for choosing  to work with this system is the different motility properties that each component exhibits. The cancer cells, represented by $c$, move using their motility apparatus on the ECM $v$ in a particular manner: they exhibit a preferred direction towards higher densities of the ECM (\textit{haptotaxis} part of the model), but still a part of their motion is random and is understood as cellular diffusion. The ECM $v$ does not move or otherwise translocate. The MMPs $m$ diffuse freely in the \textit{interstitial fluid} without possessing any motility mechanism (molecular diffusion).
				
In addition to the above motility properties, the cancer cells \textit{proliferate}  and produce MMPs with a constant rate. The MMPs attach on the ECM which they instantly dissolve (at least compared to the invasion time scale of the model), and disassemble to their constituents due to chemical degradation, see e.g.~\cite{Condeelis.2011, Rao.2003}.
				
For the needs of this paper we perform two numerical experiments in which we use the parameters
\begin{equation} \label{eq:cancer:parameters}
	\chi = 2 \times 10^{-2},~ D_c = 2 \times 10^{-4},~ D_m = 10^{-3}, ~\mu = 0.5,~\delta = 4,~ \alpha = 0.5,~ \beta = 0.3.
\end{equation}
				
For the numerical treatment we proceed in a similar way as with the Euler equations in Section \ref{sec:euler}. We solve the system~\eqref{eq:cancer:1}-\eqref{eq:cancer:2} using a FV approach, and consider a time grid $t^0=0,\,  t^{n+1}=t^n+\Delta t^n \text{ for } n=0,1,\dots$. We consider at every time $t^n$ a RSM $G^n=\left\{C_i^n,\,i=1,\dots,N^n\right\}$ where the number of computational cells $N^n$ is adapted during the computation.  The numerical solution itself is denoted at time $t^n$ by $\vec U_i^n$ for $i=1,\dots,N^n$ and is obtained through the numerical scheme
\begin{equation}\label{eq:CLScheme_II}
  \vec U_i^{n+1} = \vec U_i^n + \Delta t^n \vec S( \vec U_i^n) - \Delta t^n \sum_{C_j^n \in N(C_i^n)}
	\frac{|\partial C_{ij}^n |}{|C_i^n|} \,  \vec H(\vec U_i^n, \vec U_j^n, \vec n_{ij}^n), 
\end{equation} 
for $i=1,\dots,N^n$, where  $\partial C_{ij}^n$ denotes the edge between the cells $C_i^n$ and $C_j^n$ and $\vec{n}_{i,j}^n$ the outer normal vector of $C_i^n$ pointing towards $C_j^n$. By abuse of notation we assume in~\eqref{eq:CLScheme_II} that $ U_i^n$ already refers to the approximate solution after projection to $C_i^n$ and hence both $ U_i^n$ and $ U_i^{n+1}$ are considered on the same cell. The combined diffusion and haptotaxis flux is approximated by the numerical flux function
\begin{equation}\label{eq:cancerflux}
	\vec H(\vec U_i^n, \vec U_j^n, \pm \vec e_k) = \pm
	\begin{pmatrix}
		- D_c\, \nabla_{i,j}^h c^n & + \chi (\nabla_{i,j}^h v^n)^+ \, c_i - \chi (\nabla_{i,j}^h v^n)^- \, c_j \\ 0 &\\
		-D_m \, \nabla_{i,j}^h m^n & 
	\end{pmatrix},
\end{equation}
where $\vec e_k\in\{\vec e_1, \vec e_2\}$  refers to an unit vector in $\mathbb{R}^2$ and we use the notations
$$ \vec U_i^n = \begin{pmatrix} c_i^n \\ v_i^n \\ m_i^n
\end{pmatrix},\quad \begin{pmatrix} \nabla_{i,j}^h c^n \\ \nabla_{i,j}^h v^n \\ \nabla_{i,j}^h m^n
\end{pmatrix} = \frac{\vec U_j^n - \vec U_i^n}{|\, \partial C_{ij}^n \,|} $$
for $1 \leq i,j\leq N^n$. Moreover, the positive and negative part are defined by
$$ f^+ = \max\left\{0,~f\right\},\quad f^- = - \min\left\{0,~f\right\}.$$
\change{}{In~\eqref{eq:cancerflux} we employ classical finite differences for the discretization of diffusion, which are not considered consistent in case of hanging nodes as they appear in our AMR scheme, clarify~\cite{eymardFiniteVolumeMethods2000}. Therefore, using this discretization an additional error is introduced, which depends on the area of all cells that have a hanging node at their boundary~\cite{eymardFiniteVolumeMethods2000, belmouhoub1996}. While for this reason in general a more specialized scheme for diffusion on hanging nodes such as~\cite{coudiereConvergenceRateFinite2000} is preferable, the finite difference discretization in combination with AMR can still improve the accuracy of the full scheme when compared to a uniform method. We thus use this discretization in the present experiment for simplicity.}
				
We employ time steps $\Delta t^n$ according to the condition 
\begin{equation}\label{eq:cancercfl}
\Delta t^n \leq \text{CFL} \min_{1\leq i,j\leq N^n} \min \left\{\frac{|\partial C_{ij}^n|}{\chi |\nabla_{ij}^h v^n |}, |(\partial C_{ij}^n)|^2 \, D_m^{-1}\right\}
\end{equation}
using the CFL number $0.5$. Due to the explicit scheme we employ, the quadratic term in~\eqref{eq:cancercfl} is required for numerical stability. However, due to the low diffusivity this term does not cause a significant computational burden.  We discretize the domain using a RSM with $l_\text{min} = 5$ and $l_\text{max}=7$, start the simulation on the grid $G^0 = G_5^2$ and perform both strong refinement and weak coarsening once after each time step of scheme~\eqref{eq:CLScheme_II}. Similar as in Section~\ref{sec:euler} we consider the gradient of the cancer density 
$$g_i^n = \max_{C_j^n \in N(C_i^n)} \frac{|c_j^n - c_i^n|}{\|M(C_i^n)-M(C_j^n)\|_2}, \quad M_i^n = \frac{g_i^n}{\max_{\{ 0\leq i\leq N^n\}} g_i^n} $$
 and use the threshold values $\theta_\text{ref}=0.2,~\theta_\text{coars}= 0.1$. We project the solution between the cells as discussed in Section~\ref{sec:projections}.				
The domain $$ \Omega_\text{top} = \left\{ {\bf x} \in \Omega,~x_2 \geq \sin \left( \frac{x_1^3}{125} + \frac{2x_1 + 26}{25} + \frac{1}{20} \right) \right\},$$ determines our considered initial condition
\begin{equation} \label{eq:cancer_ini}
	\vec U_0({\bf x},t) = \begin{cases}
	(1,~0,~0.3)^T,\quad &{\bf x} \in \Omega_\text{top}\\
	(0,~v_0,~0)^T, \quad &{\bf x} \in \Omega \backslash \Omega_\text{top}
	\end{cases}.
\end{equation}
				
\begin{figure}[t]
  \pgfplotsset{colormap/YlOrRd-9}
    \begin{tikzpicture}
	\begin{groupplot}[
          /tikz/mark size=1.5pt,
          group style={group name=my plots, group size=3 by 1, horizontal sep=.8cm},
          xmin=0, xmax=1, ymin=0, ymax=1,
          ticklabel style = {font=\scriptsize},
          axis background/.style={fill=white},
          width=.35\linewidth, height=.35\linewidth, 
          ]
          \nextgroupplot[title = {$t=0$}]
          \addplot [forget plot] graphics [xmin=0, xmax=1, ymin=0, ymax=1] {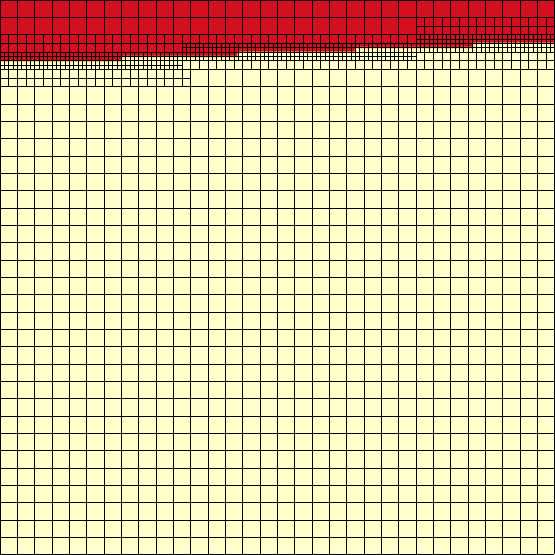};
          \nextgroupplot[title = {$t=2.5$}]
          \addplot [forget plot] graphics [xmin=0, xmax=1, ymin=0, ymax=1] {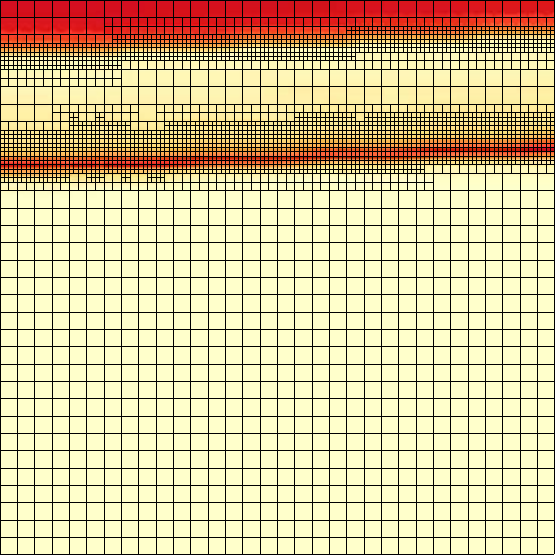};
          \nextgroupplot[title = {$t=5$},
          point meta min=0,
          point meta max=1.25,
          colorbar,
          colorbar style={at={(1.2,1)},anchor=north west}
          ]
          \addplot [forget plot] graphics [xmin=0, xmax=1, ymin=0, ymax=1] {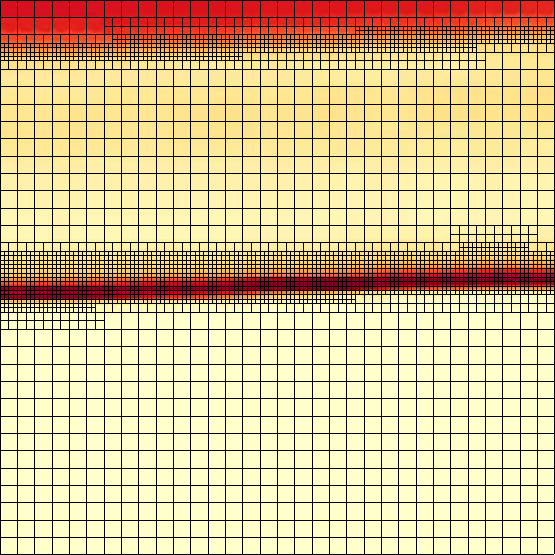};
          
         \end{groupplot}
      \end{tikzpicture}%
\caption{Numerical simulation of cancer invasion on a 2D domain with adaptive mesh refinement (time instances $t=0,$ $t=2.5$ and $t=5$ of the cancer cell density). The cancer cells invade the ECM in the form of waves and the mesh is refined at the regions of interest.}\label{fig:h.ref.cancer:1}
\end{figure}

\begin{table}[t] 
	\centering
	\begin{tabular}{l c c r}
          \toprule
		grid       & mesh cells & $L^1$ error & relative CPU time \\ \midrule
        uniform $l=5$ & $1024$ & $4.092 \times 10^{-2}$ & $4.92$ $\%$ \\
          uniform $l=6$ & $4096$ & $1.743 \times 10^{-2}$ &  $44.81$ $\%$ \\
          uniform $l=7$ & $16384$ & $7.147 \times 10^{-3}$ & $615.69$ $\%$ \\
          adaptive $l\in\{5, 6, 7\}$ & $3361$ &$1.338 \times 10^{-2}$ & $100.00$ $\% $ \\\bottomrule\\
	\end{tabular}
	\caption{$L^1$ errors and relative CPU times of the scheme in the tumor growth experiment at $t=2.5$ in case $v_0=1$. The solution computed on an uniform mesh on level $l=8$ was used as reference in the error computation. The adaptive grid solution, which uses $789$ coarse cells ($l=5$) and $396$ ($l=6$) +~$2176$ ($l=7$) refined cells, outperforms the uniform grid solutions on level $5$ and level $6$ in terms of error while taking less CPU time than the uniform grid solution on level $7$. More specifically, in order to obtain approximately half the error of the adaptive mesh solution, the uniform one on level 7 should be employed, but it implies approximately 5 times more mesh cells and approximately 6 times more CPU time.}\label{tbl:cnv:cancer:1}
      \end{table}

      \begin{table}
        \change{}{
        \begin{center}
   \begin{tabular}{l r l r}
     \toprule
     flux computation & $47.3346~\%$ & & \\ \midrule
     mesh update & $52.6654~\%$ & monitor computation & $32.4572~\%$\\ \cmidrule(l){3-4}
                      & & mesh refinement & $10.1746~\%$\\ \cmidrule(l){3-4}
                      & & mesh coarsening & $9.4059~\%$\\ \bottomrule\\
     \end {tabular}
     \end{center}
     \caption{CPU time breakdown of the AMR method in the tumor growth experiment. The computation was run in case of the uniform initial ECM density $v_0=1$ until final time $T=5$  with AMR  allowing for mesh cells from level $5$ to level $7$ (compare also Table~\ref{tbl:cnv:cancer:1}). The relative times were computed from average CPU times over 10 runs. CPU time was approximately spent half for flux computations and half for mesh updates. Within the mesh updates most CPU time was spent for the computation of the monitor function while mesh refinement and coarsening required similar CPU times. An additional test in the setting of Table~\ref{tbl:cnv:cancer:1} showed that reducing the mesh updates to one every second time step decreases the total CPU time to $73.2713~\%$ of the previous run time and the relative time for mesh updates to $35.5212~\%$, while the error was only increased by $11~\%$.}\label{tbl:breakdown}
     }
   \end{table}				
In a first tumor growth experiment, visualized in Figure \ref{fig:h.ref.cancer:1}, we consider the initial condition \eqref{eq:cancer_ini} for a constant $v_0 = 1$. We see the creation and the preliminary phase of cancer invasion in the form of waves. These cancer invasion waves emanate from the main body of the tumour and invade the ECM. They are followed by a smooth part of the solution with lower cancer cell density. Possible reconstruction of the ECM (not accounted for in the current model) would lead to a secondary wave of cancer cells invading the ECM in a similar way.  We see that the mesh is refined in the area of the first cancer wave as well as at the front of the main body of the tumour. We also note that despite the first order of accuracy of the numerical method, there is still a gain in accuracy and efficiency, by using AMR methods. This is shown in Table~\ref{tbl:cnv:cancer:1}, where the $L^1$ error (using a reference solution on a uniform grid on level 8) show that the adaptive case, where the mesh cells vary between level 5, 6 and 7, outperforms the uniform case on levels 5 and 6 at the time $t=2.5$. The $L^1$ errors have been computed using the projection to lower levels described in Section~\ref{sec:projections}. Moreover, the adaptive numerical solution is more affordable in terms of CPU time\footnote{Computations were run on a Laptop with 2.4 GHz Quad-Core Intel Core i5 CPU and 16~GB of RAM running macOS 11.4.} than the uniform one on level 7. \change{}{Moreover, a breakdown of CPU time shown in Table~\ref{tbl:breakdown} revealed that roughly half the computation time was spent for flux computations and half for the mesh updates, out of which the computation of the monitor function was the most expensive part.}
					
\begin{figure}[t] 
    \begin{tikzpicture}
	\begin{groupplot}[
          /tikz/mark size=1.5pt,
          group style={group name=my plots, group size=3 by 1, horizontal sep=.7cm},
          xmin=0.04, xmax=0.96, ymin=0.5, ymax=0.96,
          ticklabel style = {font=\scriptsize},
          xtick={0.2, 0.8},
          axis background/.style={fill=white},
          width=.4\linewidth, height=.27\linewidth, 
          ]
          \nextgroupplot[title = {ECM}]
          \addplot [forget plot] graphics [xmin=0, xmax=1, ymin=0, ymax=1] {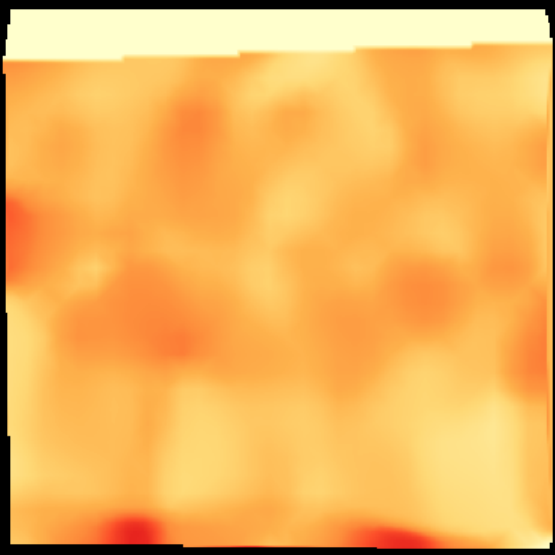};
          \nextgroupplot[title = {$t=1$}]
          \addplot [forget plot] graphics [xmin=0, xmax=1, ymin=0, ymax=1] {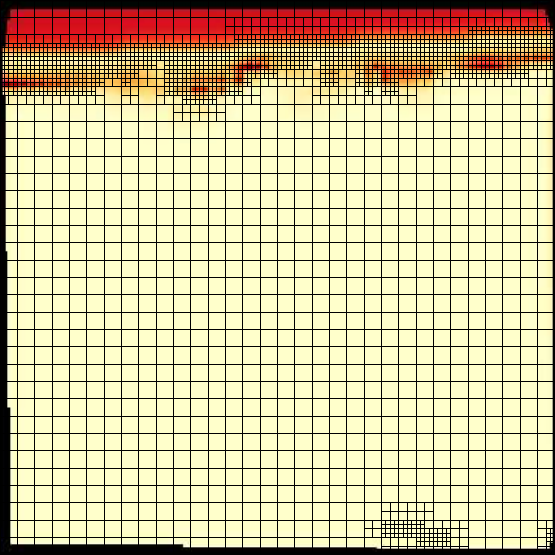};
           \nextgroupplot[title = {$t=4$},
           point meta min=0, point meta max=1.25,
          ]
          \addplot [forget plot] graphics [xmin=0, xmax=1, ymin=0, ymax=1] {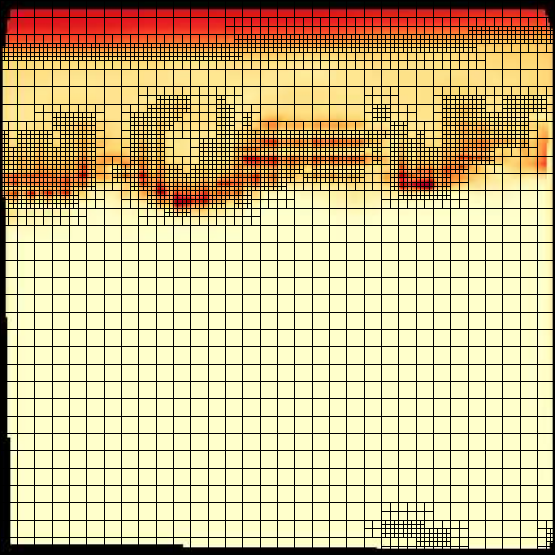};
          
         \end{groupplot}
      \end{tikzpicture}%
\caption{Numerical simulation of cancer invasion with adaptive mesh refinement on a non uniform ECM. Left panel: initial ECM density. Middle and right panel: Cancer cells density at time instances $t=1$ and $t=4$. The same colormap as in Figure~\ref{fig:h.ref.cancer:1} is used. The migration of the cancer cells is strongly influenced by the non-uniformity of the matrix. They concentrate and invade the ECM in the form of cancer cell ``islands''. The refinement of the grid follows the dynamics of the cancer cells.} \label{fig:h.ref.cancer2}
\end{figure}
					
In the second experiment, shown in Figure \ref{fig:h.ref.cancer2},  we investigate the effect that a non-uniform ECM has in the invasion of the cancer cells by imposing a spatial structure on $v_0$ in \eqref{eq:cancer_ini} shown in the left panel of Figure \ref{fig:h.ref.cancer2}. We see that, despite the simple structure of the model, the dynamics of the solution are complex. As in the previous experiment, a propagating wave is formed. This time however, due to the non-uniformity of the ECM, the cancer cells concentrate in isolated ``islands'' as they move through the matrix.  We can again notice the separation of these invading islands from the main body of the tumour; a behaviour that is consistent with the biomedical understanding of tumour spread,  \cite{JapaneseGastricCancerAssociation2011}. The mesh follows the dynamics of the solution and finely resolves the front of the main body of the tumour and the areas where the invasion takes place. The different refinement and coarsening thresholds $\theta_\text{ref}$ and $\theta_\text{coars}$ employed, cause a ``memory effect'' in the mesh, which maintains a higher resolution in previous locations of the cancer islands. These are the areas where the reconstruction of the ECM (not included in the current model) mostly takes place and a higher resolution can be very useful; another benefit of the AMR method that we employ.
					
\section{Conclusions}
				
The current work is motivated by our wider investigation of cancer invasion models. The special nature of these problems and, in particular, the highly dynamic behaviour of the solutions necessitate the development of specialized numerical methods and techniques. These methods can become expensive, especially in the multidimensional cases, and so AMR is sought as a way to decrease computational costs. The application of AMR methods in these problems, gives rise to a series of difficulties that need to be addressed.
					
We have presented here a newly developed mesh structure data administration technique used as machinery for our AMR (h-refinement) methods. When compared to existing methods in the literature, our technique exhibits similarities to pointer-based mesh data structure techniques, and exploits the rectangular structure of the mesh and its refinement by bisection.
					
We introduced an easy to use technique that avoids the traversing of the connectivity  graph for the cell ancestry and, due to the structure of the mesh, it greatly simplifies the identification of neighbouring cells. It can be easily implemented and employed in a wide range of problems in 1-, 2-, and higher dimensional spaces. It is particularly designed for smooth meshes, and uses their smoothness dynamically in the matrix operations.
The memory footprint of the method makes it affordable on coarse to very fine mesh resolutions. Additionally, our technique allows for adaptive minimum and maximum refinement levels as well as for a free choice of monitor functions and threshold parameters. Although these properties are not investigated in the current paper, they are still potentially useful in cases where the solutions exhibit multiple dynamical phenomena or blow-ups.

We have moreover presented the components of the mesh administration technique in detail and its connection to the physical discretisation of the domain. We have discussed the operations for traversing the mesh, for the identification of the sibling and neighbour cells, as well as for the local refinement and coarsening of the mesh. 
					
Finally, we have endowed this technique with an AMR method and presented its capabilities and its flexibility in three applications. The first is a generic experiment in the absence of physical or biological laws where the mesh refinement is dictated by synthetic monitor functions. The second is a physical application of the technique and the AMR in the classical case of the Euler equation. The discussion of particular numerical solvers is beyond the scope of this paper, so we have used a common vector splitting scheme. The third application is a biological problem: a 2D tumour growth and invasion of the ECM model. This model (like other cancer invasion models) exhibits highly dynamic solutions that constitute a challenge for typical numerical methods. Even for the low order of accuracy of the numerical method that we have used we have seen a gain in efficiency and accuracy obtained by the AMR technique.
					
Future steps along this direction should focus mostly on the development of higher order and---for the diffusive part of the problems---implicit numerical solvers to be used in the AMR method. The development of such methods, however, necessitates extensive analysis, algorithm implementation, and numerical experimentation, which falls beyond the scope of the current paper, and is therefore postponed for a subsequent work. 
					
\bibliographystyle{plain} 
\bibliography{new_h}
\section*{Acknowledgements}
NK was supported by the Postdoctoral Fellowships for Research in Japan (Standard) of the Japan Society for the Promotion of Science. NS was partly funded from the German Science Foundation (DFG) under the grant SFB 873:  ``Maintenance and Differentiation of Stem Cells in Development and Disease''.
\appendix

\end{document}